\documentclass[aop,MSNbibl,dvips]{arximspdf}
\usepackage{mathbh}
\usepackage{graphicx}

\doi{10.1214/10-AOP609}
\volume{39}
\issue{6}
\pubyear{2011}
\firstpage{2178}
\lastpage{2223}

\makeatletter

\let\epsilon\varepsilon

\newtheorem{defin}{Definition}
\newtheorem{theorem}{Theorem}
\newtheorem{prop}{Proposition}
\newtheorem{lemma}{Lemma}
\newtheorem{lem}{Lemma}[section]
\newproclaim{ex}{Example}
\newproclaim{rem}{Remark}
\newproclaim{claim}{Claim}

\renewcommand{\subset}{\subseteq}
\renewcommand{\supset}{\supseteq}

\newcommand{\eqref}[1]{(\ref{#1})}
\newcommand{\oa}{\overline{a}}
\newcommand{\ua}{\underline{a}}
\newcommand{\oi}{\overline{i}}
\newcommand{\ui}{\underline{i}}
\newcommand{\oip}{\overline{i+1}}
\newcommand{\uip}{\underline{i+1}}
\newcommand{\oie}{\overline{2i}}
\newcommand{\uie}{\underline{2i}}
\newcommand{\oio}{\overline{2i+1}}
\newcommand{\uio}{\underline{2i+1}}
\newcommand{\oiep}{\overline{2(i+1)}}
\newcommand{\uiep}{\underline{2(i+1)}}
\newcommand{\Di}{\operatorname{\mathsf{Diag}}}
\newcommand{\1}{\mathbh{1}}
\newcommand{\iy}{\infty}
\newcommand{\R}{\mathbb{R}}
\newcommand{\sr}{\sqrt}
\newcommand{\Z}{\mathbb{Z}}
\newcommand{\N}{\mathbb{N}}

\newcommand{\Cst}{\mathsf{Cst}}

\newcommand{\Sc}{\mathcal{S}}
\newcommand{\F}{\mathcal{F}}
\newcommand{\Hh}{\mathcal{H}}
\newcommand{\Es}{\mathbb{E}}
\newcommand{\Pb}{\mathbb{P}}
\newcommand{\Var}{\operatorname{\mathsf{Var}}}
\newcommand{\Ff}{\mathbb{F}}
\newcommand{\Gf}{\mathbb{G}}
\newcommand{\C}{\mathcal{C}}
\newcommand{\Ec}{\mathcal{E}}
\newcommand{\Pc}{\mathcal{P}}
\newcommand{\E}{\mathcal{E}}
\newcommand{\T}{\mathcal{T}}
\newcommand{\Rc}{\mathcal{R}}
\newcommand{\Nn}{\mathcal{N}}
\newcommand{\A}{\mathcal{A}}
\newcommand{\Lc}{\mathcal{L}}
\newcommand{\Gg}{\mathcal{G}}
\newcommand{\las}{\operatorname{\mathsf{Sp}}}

\newcommand{\al}{\alpha}
\newcommand{\be}{\beta}
\newcommand{\G}{\Gamma}
\newcommand{\De}{\Delta}
\newcommand{\Ep}{\Upsilon}
\newcommand{\ze}{\zeta}
\newcommand{\Tet}{\Theta}
\newcommand{\ki}{\chi}
\newcommand{\ka}{\kappa}
\newcommand{\La}{\Lambda}
\newcommand{\s}{\sigma}
\newcommand{\Om}{\Omega}

\def\des{\longrightarrow}
\def\ses{\sim}

\makeatother

\begin{document}
\begin{frontmatter}

\title{Dynamics of vertex-reinforced random walks}
\runtitle{Dynamics of vertex-reinforced random walks}

\begin{aug}
\author[A]{\fnms{Michel} \snm{Bena{\"{\i}}m}\ead[label=e1]{michel.benaim@unine.ch}\thanksref{aut1}}
and
\author[B]{\fnms{Pierre} \snm{Tarr\`{e}s}\corref{}\ead[label=e2]{tarres@math.univ-toulouse.fr}\thanksref{aut2}}
\thankstext{aut1}{Supported by the Swiss National Foundation Grant
200020-120218/1.}
\thankstext{aut2}{On leave from the University of Oxford.
Supported in part by the Swiss National Foundation Grant
200021-1036251/1, and by a Leverhulme Prize.}
\runauthor{M. Bena{\"{\i}}m and P. Tarr\`{e}s}
\affiliation{Universit\'{e} de Neuch\^{a}tel and CNRS, Universit\'{e}
de Toulouse}
\address[A]{Institut de math\'{e}matiques\\
Universit\'{e} de Neuch\^{a}tel\\
Rue Emile-Argand 11, Case postale 2\\
2007 Neuch\^{a}tel\\
Switzerland\\
\printead{e1}} 
\address[B]{Institut de Math\'{e}matiques\\
CNRS, Universit\'{e} de Toulouse\\
118 route de Narbonne\\
31062 Toulouse Cedex 9\\
France\\
\printead{e2}}
\end{aug}

\received{\smonth{11} \syear{2008}}
\revised{\smonth{9} \syear{2010}}

%
\begin{abstract}
We generalize a result from Volkov [\textit{Ann. Probab.} \textbf{29}
(2001) \mbox{66--91}] and prove that,
on a large class of locally finite connected graphs of bounded degree
$(G,\sim)$ and symmetric reinforcement matrices $a=(a_{i,j})_{i,j\in
G}$, the vertex-reinforced random walk (VRRW) eventually localizes with
positive probability on subsets which consist of a complete $d$-partite
subgraph with possible loops plus its outer boundary.

We first show that, in general, any stable equilibrium of a linear
symmetric \textit{replicator }dynamics with positive payoffs on a graph
$G$ satisfies the property that its support is a complete $d$-partite
subgraph of $G$ with possible loops, for some $d\ge1$. This result is
used here for the study of VRRWs, but also applies to other contexts
such as evolutionary models in population genetics and game theory.

Next we generalize the result of Pemantle [\textit{Probab. Theory
\mbox{Related} Fields} \textbf{92} (1992) 117--136] and
Bena{\"{\i}}m [\textit{Ann. Probab.} \textbf{25} (1997) 361--392] relating the asymptotic behavior of
the VRRW to \textit{replicator } dynamics. This enables us to conclude
that, given any neighborhood of a strictly stable equilibrium with
support $S$, the following event occurs with positive probability: the
walk localizes on $S\cup\partial S$ (whe\-re~$\partial S$ is the outer
boundary of $S$) and the density of occupation of the VRRW converges,
with polynomial rate, to a strictly stable equilibrium in this
neighborhood.
\end{abstract}

%

\end{frontmatter}

\section{General introduction}
\label{sec:gintro}
Let $(\Om,\F,P)$ be a probability space.
Let $(G,\sim)$ be a locally finite connected symmetric graph, and let
$G$ be its vertex set, by a slight abuse of notation.
Let $a:=(a_{i,j})_{i,j\in G}$ be a symmetric (i.e., $a_{i,j}=a_{j,i}$)
matrix with nonnegative entries
such that, for all $i, j \in G$,
\[
i\sim j \quad \Leftrightarrow\quad  a_{i,j}>0.
\]
Let $(X_n)_{n\in\N}$ be a process taking values in $G$.
Let $\Ff=(\F_n)_{n\in\N}$ denote the filtration
generated by the process, that is,
$\F_n=\s(X_0,\ldots,X_n)$ for all $n\in\N$.

For any $i\in G$, let $Z_n(i)$ be the number of times
that the process visits site $i$ up through
time $n\in\N\cup\{\iy\}$, that is,
\[
Z_n(i)=Z_0(i)+\sum_{m=0}^n\1_{\{X_m=i\}}
\]
with the convention that, before initial time $0$, a site $i\in G$ has
already been visited $Z_0(i)\in\R_+\setminus\{0\}$ times.

Then $(X_n)_{n\in\N}$ is called a \textit{Vertex-Reinforced Random
Walk}
(\textit{VRRW})
\textit{with starting point $v_0\in G$ and reinforcement matrix}
$a:=(a_{i,j})_{i,j\in G}$
if
$X_0=v_0$ and, for all $n\in\N$,
\[
\Pb(X_{n+1}=j | \F_n)=\1_{\{j\sim X_n\}}\frac{a_{X_n,j}Z_n(j)}
{\sum_{k\sim X_n} a_{X_n,k}Z_n(k)}.
\]

These non-Markovian random walks
were introduced in 1988 by Pemantle \cite{pemantle3} during his PhD
with Diaconis, in the spirit of the model of
Edge-Reinforced Random Walks by Coppersmith and Diaconis in 1987 \cite
{diaconis}, where the weights accumulate on edges rather than vertices.

Vertex-reinforced random walks were first studied in the articles of
Pemantle \cite{pemantle1} and Bena{\"{\i}}m \cite{benaim1} exploring
some features of their asymptotic behavior on
finite graphs and, in particular, relating the behavior of the
empirical occupation measure to solutions of ordinary differential
equations when the graph is complete (i.e., when all vertices are
related together), as explained below. On the integers $\Z$, Pemantle
and Volkov
\cite{pemantle2} showed that the VRRW a.s. visits only finitely
many vertices and, with positive probability, eventually gets stuck on
five vertices, and Tarr\`{e}s \cite{tarres} proved that this
localization on five points is the almost sure behavior.

On arbitrary graphs, Volkov \cite{volkov} proved that VRRW with
reinforcement coefficients $a_{i,j}=\1_{i\sim j}$, $i, j\in G$
(again, $i\sim j$ meaning that $i$ and $j$ are neighbors in the
nonoriented graph $G$), localizes with positive probability on some
specific finite subgraphs; we recall this result in Theorem \ref
{thm:traps} below, in a generalized version. More recently, Limic and
Volkov \cite{limic-volkov} study VRRW with the same specific type of
reinforcement on complete-like graphs (i.e., complete graphs ornamented
by finitely many leaves at each vertex) and show that, almost surely,
the VRRW spends positive (and equal) proportions of time on each of its
nonleaf vertices.

The VRRW with polynomial reinforcement [i.e., with the probability to
visit a vertex proportional to a function $W(n)=n^\rho$ of its current
number of visits] has recently been studied by Volkov on $\Z$
\cite{volkov2}. In the superlinear
case (i.e., $\rho> 1$), the walk a.s. visits two vertices infinitely
often. In the sublinear case (i.e., $\rho< 1$), the walk a.s. either visits
infinitely many sites infinitely often or is transient; it is
conjectured that the latter behavior cannot occur, and that, in
fact, all integers are infinitely\vadjust{\goodbreak} often visited.

The similar Edge-Reinforced Random Walks and, more generally,
self-interact\-ing processes, whether in discrete or continuous
time/space, have been extensively studied in recent years. They are
sometimes used as models involving
self-organization or learning behavior, in physics, biology or
economics. We propose a~short review of the subject in the
introduction of~\cite{mountford-tarres}. For more detailed overviews,
we refer the reader to surveys by Davis~\cite{davis3}, Merkl and
Rolles~\cite{merkl-rolles3}, Pemantle~\cite{pemantle5} and T\'{o}th~\cite{toth},
each analyzing the subject from a different perspective.

Let us first recall a few well-known observations on the study of
Vertex-Reinforced Random Walks, and, in particular,
the heuristics for relating its behavior to solutions of ordinary
differential equations when the graph is finite and complete (i.e.,
when all
vertices are related together), as done in Pemantle \cite{pemantle1}
and Bena{\"{\i}}m \cite{benaim1}.

Let us introduce some preliminary notation, without any further  assumption on
$(G,\sim)$ locally finite connected symmetric graph, possibly
infinite. For all $x=(x_i)_{i\in G}\in\R^{G}$, let
\[
S(x):=\{i\in G/ x_i\not=0\}
\]
be its support. For all $x\in\R^{G}$ such that $S(x)$ is finite, let
%
\begin{equation}
\label{nih}
N_i(x):=\sum_{j\in G, j\sim i}a_{i,j}x_j,
H(x)=\sum_{i,j\in G, i\sim j}a_{i,j}x_ix_j=\sum_{i\in G}x_iN_i(x)
\end{equation}
and, if $H(x)\not=0$, let
%
\begin{equation}
\label{pi}
\pi(x):= \biggl(\frac{x_iN_i(x)}{H(x)} \biggr)_{i\in G}.
\end{equation}
Let
\[
\Tet:= \{x\in\R^{G}\mbox{ s.t. } |S(x)|<\iy\},
\]
and let
\[
\De:= \biggl\{x\in\R_+^{G}\cap\Tet\mbox{ s.t. }\sum_{i\in G}x_i=1 \biggr\}
\]
be the nonnegative simplex restricted to elements $x$ of finite support.

For all $n\in\N$, let
\[
y(n):=\biggl(\frac{Z_n(i)-Z_0(i)}{n}\biggr)_{i\in G}\in\Tet\cap\De
\]
and, if $G$ is finite, let
\[
x(n):=\biggl(\frac{Z_n(i)}{n+n_0}\biggr)_{i\in G},
\]
where $n_0:=\sum_{j\in G}Z_0(j)>0$: $x(n)$ [resp. $y(n)$] is the
vector of density of occupation  of the random walk at time $n$, with
the convention that site $i$ has been visited $Z_0(i)$ (resp., $0$)
times at\vadjust{\goodbreak} time $0$.

Assume, for the
sake of simplicity in the following heuristic argument, that~$G$ is a
finite graph. Let $L\gg1$. For all $n\in\N$, the goal is to compare $x(n+L)$ to $x(n)$.
If $n\gg L$, then the VRRW between these times behaves as though~$x(k)$, $n\le k\le n+L$,
were constant, and hence approximates a~Markov chain which we
call $M(x(n))$.

Then $\pi(x(n))\in\De$ is the invariant measure of
$M(x(n))$, which is reversible [trivially $H(x(n))>0$ since $x(n)_i>0$
for all $i$, so that $\pi(x(n))$ is well defined].
If $L$ is large enough, then, by the ergodic theorem, the local
occupation density between these times will be close to $\pi(x(n))$.
This means that
%
\begin{equation}
\label{mchain}
(n+L)x(n+L)\approx nx(n)+L \pi(x(n)),
\end{equation}
hence,
%
\begin{equation}
\label{mdyn}
x(n+L)-x(n)\approx\frac{L}{nH(x(n))}F(x(n)),
\end{equation}
where
%
\begin{equation}
\label{replicator}
F(x)=(x_i[N_i(x)-H(x)])_{i\in G}.
\end{equation}

Up to an adequate time change, $(x(k))_{k\in\N}$ should approximate
solutions of the ordinary differential equation on $\De$,
%
\begin{equation}
\label{eqdiff}
\frac{dx}{dt}=F(x),
\end{equation}
also known as the linear \textit{replicator }equation in population
genetics and game theory.

However, the requirement that $L$ be large enough so that the local
occupation measure of the Markov Chain approximates the invariant
measure~$\pi(x(n))$ competes with the other requirement that $L$ be
small enough so that the probability transitions of this Markov Chain
still match the ones of the VRRW, so that the heuristics breaks down
when the relaxation time of the Markov Chain is of the order of $n$,
which can happen in general on noncomplete graphs and is actually
consistent with the fact that the walk will indeed eventually localize
on a small subset. An illustration of how such a behavior can occur is
given in the proof of Lemma 2.8 in Tarr\`{e}s \cite{tarres}. The study
of the a.s. asymptotic behavior of the VRRW on an infinite graph is
even more involved in general.

Let us yet study the replicator differential equation (\ref{eqdiff})
associated to the random walk on $\De$ for general locally finite
symmetric graphs $(G,\sim)$.

It is easy to check that $H$ is a strict Lyapounov function for (\ref
{eqdiff}) on $\De$, that is,  strictly increasing on the
nonconstant solutions of this equation: if $x(t)=(x_i(t))_{i\in G}$ is
the solution at time $t$, starting at $x(0):=x_0$, then
\[
\frac{d}{dt}H(x(t))=
\sum_{i\in S(x)}\frac{\partial H}{\partial x_i}(x(t))F(x(t))_i=J(x(t)),
\]
where, for all $x\in\De$,
%
\begin{equation}
\label{defj}
J(x):=2\sum_{i\in S(x)}N_i(x)F(x)_i=2\sum_{i\in S(x)} x_i\bigl(N_i(x)-H(x)\bigr)^2.
\end{equation}
Note that the restriction of $H$ to the equilibria of (\ref{eqdiff})
takes finitely many values if $G$ is finite (see \cite{pemantle1}, e.g.).\vadjust{\goodbreak}

Let us now deal with the equilibria of this differential equation: a point
$x=(x_i)_{i\in G}\in\De$ is called an \textit{equilibrium} if and
only if
$F(x) = 0.$ An equilibrium is called \textit{feasible} provided
$H(x)\not=0.$

On a finite graph $G$, any equilibrium point $x\in\De$ of
$(x(n))_{n\in\N}$ is feasible: for all $n\in\N$ and $i\in G$,
$Z_n(i)\le\sum_{j\sim i}Z_n(j)+n_0$, so that $x$ would satisfy $N_i(x)\ge(\min_{j\sim i}a_{i,j})x_i$ for all
$i\in G$, hence,
%
\begin{equation}
\label{lbound_H}
H(x)\ge\Bigl(\min_{\{i,j\in S(x), j\sim i\}}a_{i,j} \Bigr)\sum_{i\in
S(x)}x_i^2\ge\frac{\min_{\{i,j\in S(x), j\sim i\}}a_{i,j}}{|S(x)|}
\end{equation}
by the Cauchy--Schwarz inequality.

By a slight abuse of notation, we let $DF(x) = (\partial F_i/\partial
x_j)_{i,j\in G}$ denote both the \textit{Jacobian matrix} of $F$ at $x$,
and the corresponding linear operator on $\Tet$. Since $\Delta$ is
invariant under the flow induced by $F,$ the tangent space
\[
T\De:= \biggl\{y\in\Tet\Big/ \sum_{i\in G}y_i=0 \biggr\}
\]
is invariant under $DF(x).$ We let $DF(x)|_{T\Delta}$ denote the
restriction of the operator $DF(x)$ to
$T\Delta.$

When $x$ is an equilibrium, it is easily seen that $DF(x)$ has real
eigenva\-lues (see Lemma \ref{lem:1}). Such an equilibrium is called
\textit{hyperbolic} (resp., \textit{a sink}) provided
$DF(x)|_{T\Delta
}$ has nonzero (resp., negative) eigenvalues. It is called a~\textit
{stable equilibrium} if $DF(x)|_{T\De}$ has nonpositive eigenvalues.
%
Note that every sink is stable. Furthermore, by Theorem \ref
{partition} below, every stable equilibrium is feasible.

We will sometimes abuse notation and identify arbitrary subsets $H$ of
$G$ to the corresponding subgraph $(H,\sim)$.
Given $i\in G$ and a subset $A$ of $G$, we write $i\sim A$ if there
exists $j\in A$ such that $i\sim j$. Given two subsets $R$ and $S$ of
$G$, we let
\[
\partial R=\{j\in G\setminus R \dvtx j\sim R\},\qquad
\partial_{S} R=\{j\in S\setminus R \dvtx j\sim R\};
\]
$\partial R$ is called the \textit{outer boundary} of $R$.

Given $e$, $e'\in E(G)$, we write $e\sim e'$ if $e$ and $e'$ have
at least one vertex in common.

A site $i\in G$ will be called a loop if $i\sim i$, and we will say
that a subset $H$ contains a loop iff there exists a site in it which
is a loop.

We will say that $x$ is a \textit{strictly stable equilibrium} if it is
stable and, furthermore, for all $i\in\partial S(x)$, $N_i(x)<H(x).$
We let $\Ec_s$ be the set of strictly stable equilibria of (\ref
{eqdiff}) in $\De$. Note that $x$ stable already implies $N_i(x)\le
H(x)$ for all $i\in\partial S(x)$, by Lemma \ref{lem:1}.

Given $d\ge1$, subgraph $(S,\sim)$ of $(G,\sim)$ will be called a
\textit{complete $d$-partite graph with possible loops}, if $(S,\sim
)$ is
a $d$-partite graph on which some loops have possibly been added. That is,
\[
S=V_1\cup\cdots\cup V_d
\]
with:
\begin{longlist}[(ii)]
\item[(i)]
$\forall p\in\{1,\ldots,d\}$, $\forall i, j \in V_p$, if
$i\not=j$ then $i\not\sim j.$\vadjust{\goodbreak}

\item[(ii)] $\forall p, q \in\{1,\ldots,d\}$, $p\not=q$,
$\forall
i\in V_p$, $\forall j\in V_q$, $i\sim j$.
\end{longlist}

For all $S\subset G$, let (P)$_S$ be the following predicate:
\begin{enumerate}[(P)$_S\mathrm{(a)}$]
\item[(P)$_S\mathrm{(a)}$]
$(S,\sim)$ is a complete $d$-partite graph with possible loops.
\item[(P)$_S\mathrm{(b)}$]
If $i\sim i$ for some $i\in S$, then the partition containing $i$ is a
singleton.
\item[(P)$_S\mathrm{(c)}$]
If $V_p$, $1\le p\le d$ are its $d$ partitions, then for all
$p, q \in\{1,\ldots,d\}$ and
$i, i' \in V_p$, $j, j' \in V_q$,
$a_{i,j}=a_{i',j'}$.
\end{enumerate}

In the following Theorems \ref{partition}--\ref{thm:traps} and
Propositions \ref{lem:str1} and \ref{prop:str2}, we only assume the
graph $(G,\sim)$ to be symmetric and locally finite, without any
further conditions than the ones mentioned in the statements.

\begin{theorem}
\label{partition}
If $x\in\De$ is a stable equilibrium of (\ref{eqdiff}), then $x$ is
feasible and \textup{(P)}$_{S(x)}$ holds.
\end{theorem}

In the case $a=(a_{i,j})_{i,j\in G}=(\1_{i\sim j})_{i,j\in G}$ the
following Theorem \ref{partition2} provides a necessary and sufficient
condition for $x\in\De$ being a stable equilibrium. Theorems~\ref
{partition} and \ref{partition2} are proved in Section \ref{sec:spartition}.

\begin{theorem}
\label{partition2}
Assume $a_{i,j}=\1_{i\sim j}$ for all $i, j \in G$, and let
$x=(x_i)_{i\in G}\in\De$.

If $(S(x),\sim)$ contains no loop, then $x$ is a stable (resp.,
strictly stable) equilibrium if and only if there exists $d\ge2$ such
that:
\begin{longlist}[(iii)]
\item[(i)]$\!\!(S(x),\!\sim)$ is a complete $d$-partite subgraph, with $\mbox{partitions}\,{=:}\,V_1, \ldots,V_d$,
\item[(ii)]$\!\!\sum_{i\in V_p} x_i=1/d$ for all $p\in\{1,\ldots,d\}$,\vspace*{2pt}
\item[(iii)]$\!\!\forall i\in\partial S(x), N_i(x)\le$ (resp., $<$) $1-1/d$.
\end{longlist}

If $(S(x),\sim)$ contains a loop, then $x$ is a stable (resp., strictly
stable) equilibrium if and only if $(S(x),\sim)$ is a clique of loops
[resp., with the additional assumption: $\forall j\in\partial S(x)$,
$N_j(x)<1$ or, equivalently, $\partial\{j\}\not\supseteq S(x)$].
\end{theorem}
\begin{rem}
\label{rem:1}
Jordan \cite{jordan} independently shows, in the context of
preferential duplication graphs, that conditions (i)--(iii) in
Theorem \ref{partition2} are indeed sufficient for $x\in\De$ being a
stable equilibrium when loops are not allowed.
\end{rem}
\begin{rem}
\label{rem:2}
A connection between the number of stable rest points in the
\textit{replicator} dynamics [or of patterns of evolutionary stable sets
(ESS's)] and the numbers of cliques of its graph was made by Vickers
and Cannings~\cite{vickers1,vickers2}, Broom et al. \cite{broom} and
Tyrer et al. \cite{cannings}, motivated by the study of evolutionary
dynamics in biology.

A consequence of Theorem \ref{partition} is that supports of stable
equilibria are \textit{generically} cliques of the graph $G.$ More
precisely, assume that the coeffi\-cients~$(a_{i,j})_{i,j\in G}$ are
distributed according to some absolutely continuous distribution w.r.t.
the Lebesgue measure on symmetric matrices. Then the supports of stable
equilibria are a.s. cliques of the graph $G$ (i.e., any two different
vertices are connected), as a consequence of (P)$_{S(x)}$(a)
and (c).
\end{rem}

The following Theorem \ref{thm:conv} states that, given any
neighborhood $\Nn(x)$ of a~strictly stable equilibrium $x\in\Ec_s$, then,
with positive probability, the VRRW eventually localizes in
\[
T(x):=S(x)\cup\partial S(x),
\]
and the vector of density of occupation converges toward a point in
$\Nn(x)$, which will not necessarily be $x$ (there may exist a
submanifold of stable equilibria in the neighborhood of $x$).
Note that this will imply, using Remark \ref{rem:2}, that the VRRW
generically localizes with positive probability on subgraphs which
consist of a clique plus its outer boundary.

More precisely, let us first introduce the following definitions.
For all \mbox{$R\subset G$}, let
\[
\Sc(R):=S^{-1}(R)\cap\De=\{x\in\De\mbox{ s.t. }S(x)=R\}.
\]
For any open subset $U$ of $\De$ containing $x\in\De$, let $\Lc(U)$
be the event
\begin{eqnarray*}
&&\Lc(U):=\Bigl\{y(\iy):=\lim_{n\to\iy}y(n)\mbox{ exists (coordinatewise) and belongs to}\\
&&\hspace*{226pt} \Ec_s\cap\Sc(S(x))\cap U\Bigr\}.
\end{eqnarray*}

Let $\Rc$ be the asymptotic range of the VRRW, that is,
\[
\Rc:=\{i\in G\mbox{ s.t. }Z_\iy(i)=\iy\}.
\]

For any random variable $x$ taking values in $\De$, let
\[
\A_\partial(x):= \biggl\{\!\forall i\in\partial S(x), \frac
{Z_n(i)}{n^{N_i(x)/H(x)}}\mbox{ converges to a (random) limit }\in
(0,\iy) \!\biggr\}.
\]

\begin{theorem}
\label{thm:conv}
Let $x\in\De$ be a strictly stable equilibrium. Then, for any open
subset $U$ of $\De$ containing $x$,
\[
\Pb\bigl(\{\Rc=T(x)\}\cap\Lc(U)\cap\A_\partial(y(\iy))\bigr)>0.
\]
Moreover, the rate of convergence is at least reciprocally polynomial,
that is, by possibly restricting the neighborhood $U$ of $x$, there
exists $\nu :=\Cst(x,a)$ such that, a.s. on $\Lc(U)$,
\[
\lim_{n\to\iy}\bigl(y(n)-y(\iy)\bigr)n^{\nu}=0.
\]
\end{theorem}

Theorem \ref{thm:conv} is proved in Section \ref{sec:thm-conv}. It
naturally leads to the following questions.

First, are all the trapping subsets always of the form $T(x)$ for some~\mbox{$x\in\E_s$}?
The answer is negative in general: let us consider, for
instance, the graph \mbox{$(\Z,\sim)$} of integers, to which we add a loop
$0\sim0$ at site $0$, with $a_{i,j}:=\1_{i\sim j}$. Then $x:=(\1_{\{
i=0\}})_{i\in\Z}$ is a stable equilibrium, but is not strictly stable
since $N_{-1}(x)=N_1(x)=1=H(x)$. However, Proposition \ref{prop:cex} (proved
in Appendix~\ref{sec:cex}) shows that $y(n)$ converges to $x$ with
positive probability, by combining an urn result from Athreya \cite
{athreya}, Pemantle and Volkov \cite{pemantle2} (Theorem 2.3) with
martingale techniques from Tarr\`{e}s \cite{tarres} (Section 3.1).

\begin{prop}
\label{prop:cex}
Let $(G,\sim)$ be the graph of integers defined above, and~let $a_{i,j}:=\1
_{i\sim j}$. Then, with positive probability, the VRRW localizes on $\{
-2,-1,\allowbreak 0, 1,2\}$, and there exist random variables $\al\in(0,1)$, $C$
and $C' >0$ such~that
\begin{eqnarray*}
\mbox{\textup{(i)}} &\hspace*{-4pt}\quad&\frac{Z_n(0)}{n}\des_{n\to\iy}1,\\
\mbox{\textup{(ii)}} &\hspace*{-4pt}\quad&\frac{(Z_n(-1),Z_n(1))}{n/\log n}\des
_{n\to\iy}(\al,1-\al
),\\
\mbox{\textup{(iii)}}&\hspace*{-4pt}\quad& \biggl(\frac{Z_n(-2)}{(\log n)^\al
},\frac{Z_n(2)}{(\log n)^{1-\al
}} \biggr)\des_{n\to\iy}(C,C').
\end{eqnarray*}
\end{prop}
We conjecture that, conditionally on a localization of the VRRW on a
finite subset, its vector of density of occupation on the subset
converges to a~stable equilibrium $x$ of (\ref{eqdiff}), that the
asymptotic range $\Rc$ is a subset of $S(x)\cup\partial S(x)\cup
\partial(\partial S(x))$, and is equal to $T(x)=S(x)\cup\partial
S(x)$ if $x\in\E_s$, which occurs generically on $a$ (in the sense
given in the paragraph after Remark~\ref{rem:2}).

A proof would require a deeper understanding of the dynamics of
$(Z_.(i))_{i\in G}$ (see Lemma \ref{bound_rn}). Note that, on the
integers $\Z$ with standard adjacency---unli\-ke Proposition \ref
{prop:cex}---and with $a_{i,j}=\1_{i\sim j}$, the result that the VRRW
a.s.~locali\-zes on five sites \cite{tarres} implies that only
equilibria in $\E_s$ are reached with positi\-ve probability. More
precisely, in this case there exist a.s. $k\in\Z$ and $x\in\De$
with $x_k=1/2$, $x_{k-1}=\al/2$, $x_{k+1}=(1-\al)/2$, $\al\in(0,1)$
(thus, $x\in\E_s$) such that $Z_n(i)/n\,{\to}\,x_i$ as $n\,{\to}\,\iy$ for all
$i\,{\in}\,\Z$, $\A_\partial(x)$ holds and $\Rc\,{=}\,T(x)$; see~\cite
{tarres}. Stable equilibria which are not in $\E_s$ correspond to
cases $\al=0$ or~$1$, which would lead to localization on six vertices
if they were possible, similarly to Proposition \ref{prop:cex}. This
result on $\Z$ can be related to the property that every neighborhood
of any stable equilibrium $x$ contains a strictly stable one.\vadjust{\goodbreak}

Second, which subsets are of the form $T(x)=S(x)\cup\partial S(x)$ for
some~\mbox{$x\in\E_s$}? We know from Theorem \ref{partition} that subsets
$S(x)$ satisfy (P)$_{S(x)}$ and thus always consist of a complete
$d$-partite subgraph with possible loops and its outer boundary for
some $d\ge2$. But (P)$_{S(x)}$ is not sufficient, and the
occurrence of such subsets also depends on the reinforcement matrix
$a=(a_{i,j})_{i,j\in G}$. Even in the case $a=(a_{i,j})_{i,j\in G}=(\1
_{i\sim j})_{i,j\in G}$ Theorem \ref{partition2} provides explicit
criteria for $x\in\E_s$, but the corresponding condition (iii)
[when $(S(x),\sim)$ has no loops] is on $x$, thus not explicitly on
the subgraph.

We introduce in the following Definition \ref{def:traps} the notion of
\textit{strongly trapping subsets}, which we prove in Theorem \ref
{thm:traps} to always be such subsets $T(x)$ for some $x\in\E_s$.
As a consequence, by Theorem
\ref{thm:conv}, the VRRW localizes on these subsets with positive
probability. The result is thus a generalization to arbitrary
reinforcement matrices of Theorem 1.1 by Volkov \cite{volkov}
when $a_{i,j}:=\1_{\{i\sim j\}}$, in which case the assumptions of
Definition \ref{def:traps} obviously reduce to (c) or (c)$'$.

\begin{defin}
\label{def:traps}
A subset $T\subset G$ is called a strongly trapping subset of $(G,\sim
)$ if $T=S\cup\partial S$, where:
\begin{longlist}[(a)]
\item[(a)] $(i,j)\mapsto a_{i,j}$ is constant on $\{(i,j)\in S^2$
s.t. $i\sim j\}$, with common\break \mbox{$value =:a_S$},
\item[(b)] $\max_{i\in S, j\in\partial S}a_{i,j}\le a_S$, and
\end{longlist}
either
\begin{longlist}[(c)(ii)]
\item[(c)(i)]
$S$ is a complete $d$-partite subgraph of $G$ for some $d\ge2$,
 with partitions $V_1, \ldots, V_d$,
\item[(ii)] $\forall j\in\partial S$,
$\exists p\in\{1,\ldots,d\}$ and $i\in S\setminus V_p$ such that
$j\not\sim V_p\cup\{i\}$,
\end{longlist}
or
\begin{longlist}[(c)$'$]
\item[(c)$'$] $S$ is a clique of
loops, and $\forall j \in\partial S$, $\partial\{j\}\not\supset S$.
\end{longlist}
\end{defin}
\begin{theorem}
\label{thm:traps}
Let $T$ be a strongly trapping subset of $(G,\sim)$; then the VRRW has
asymptotic range $T$ with positive probability.

More precisely, assume $T=S\cup\partial S$, where $S$ satisfies
conditions \textup{(a)--(c)} or \textup{(c)}$'$ of Definition \ref{def:traps}, and
let us use the corresponding notation. Let
\begin{eqnarray*}
\Sigma&:=& \biggl\{x\in\Sc(S)\mbox{ s.t. }\sum_{i\in V_q}x_i=1/d\mbox{
for all }1\le q\le d \biggr\},\\
r_d &:=&d/(d-1)
\end{eqnarray*}
if $(S,\sim)$ contains no loops, and $\Sigma:=\Sc(S)$, $r_d:=1$ otherwise.

\begin{figure}

\includegraphics{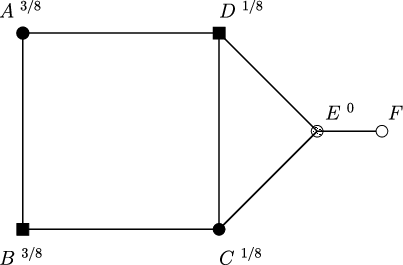}

\caption{We show in Example
\protect\ref{ex:1} that $T:=\{A,B,C,D,E\}$ does not satisfy the
assumptions of Theorem \protect\ref{thm:traps}, but is a trapping
subgraph with
positive probability by Theorems \protect\ref{partition2} and
\protect\ref{thm:conv}.
The numbers indicated in superscript of vertices represent the limit
proportions of visits to these vertices if $x(n)$ were to converge to
the equilibrium $x$ in the example. In this case the walk would
asymptotically spend most of the time in the bipartite subgraph
$S:=V_1\cup V_2$, where $V_1:=\{A,C\}$, $V_2:=\{B,D\}$, evenly divided
between partitions $V_1$ and $V_2$, and vertex~$E$ would be seldom
visited, of the order of $\sqrt{n}$ times at time $n$.}
\label{fig1}\label{counterexample}
\end{figure}

Then, for any $x\in\Sigma$ and any neighborhood $\Nn(x)$ of $x$ in
$\Sigma$, there exist random variables $y\in\Nn(x)$ and $C_j>0$,
$j\in\partial S$ such that, with positive probability:
\begin{longlist}[(iii)]
\item[(i)] VRRW eventually localizes on $T$, that
is,\vadjust{\goodbreak} $\Rc=T$,
\item[(ii)] $Z_n(i)/n\des_{n\to\iy} y_i$ for all
$i\in S$,
\item[(iii)] $Z_n(j)\ses_{n\to\iy} C_jn^{r_d\sum
_{i\sim j}a_{i,j}y_i/a_S}$ for all $j\in\partial S$.
\end{longlist}
\end{theorem}

Theorem \ref{thm:traps} is proved in Section \ref{sec:thm-traps}. We
provide in Example \ref{ex:1} (illustrated in Figure \ref{counterexample}) a
counterexample showing that Theorem \ref{thm:conv} is stronger, even
in the case $a=(\1_{i\sim j})_{i,j\in G}$.

Third, which conditions on the graph and on the reinforcement matrix~$a$
do ensure the existence of at least one strictly stable equilibrium
$x\in\E_s$, thus implying localization with positive probability on
$T(x)$? First note that, trivially, this does not always occur, for
instance, on $\Z$ when $\phi(n):=a_{\{n,n+1\}}$ is strictly monotone,
in which case we believe the walk to be transient.

In the case $a=(\1_{i\sim j})_{i,j\in G}$, Volkov \cite{volkov}
proposed the following result, using an iterative construction on
subsets of the graph.
\begin{prop}[(Volkov \cite{volkov})]
\label{lem:str1}
Assume that $a=(\1_{i\sim j})_{i,j\in G}$, and that $(G,\sim)$ does
not contain loops. Then, under either of the following conditions,
there exists at least one strongly trapping subset:
\begin{longlist}[(C)]
\item[(A)] $(G,\sim)$ does not contain triangles;
\item[(B)] $(G,\sim)$ is of bounded degree;
\item[(C)] the size of any complete subgraph is uniformly bounded by
some number~$K$.
\end{longlist}
\end{prop}
\begin{pf}
Start, for some $d\ge2$, with any complete $d$-partite subgraph
$(S,\sim)$ of $G$ with partitions $V_1, \ldots, V_d$ (e.g.,
a pair of connected vertices, $d=2$). Let $x\in\partial S$,
$S=V_1\cup\cdots\cup V_d$:
\begin{longlist}[(3)]
\item[(1)] First assume that $x\sim V_p$ for all $1\le p\le d$. Then,
for all $1\le p\le d$, let $j_p\in V_p$ be such that $x\sim j_p$;
iterate the procedure with the subgraph $\bigcup_{1\le p\le
d}\{j_p\}\cup\{x\}$, which is a clique, and thus a complete
$(d+1)$-partite subgraph.

\item[(2)] Now assume there exists $p$ such that $x\not\sim V_p$, with
$\partial\{x\}\supset S\setminus V_p$. Then we iterate the procedure
with the complete $d$-partite subgraph $S\cup\{x\}$ with partitions
$V_1, \ldots, V_p\cup\{x\}, \ldots, V_d$.

\item[(3)] Otherwise we keep the same subgraph $S$ and try another
$x\in\partial S$.
\end{longlist}

The construction eventually stops if (A), (B) or (C)
holds. When it does, that is, when $S$ has remained unchanged for all
$x\in\partial S$, then $T=S\cup\partial S$ is a strongly trapping
subgraph in the sense of Definition \ref{def:traps}.
\end{pf}

Using a similar technique, we can obtain the following necessary
condition for the existence of a strongly trapping subset in the case
of general reinforcement matrices $a$, when the graph does not contain
triangles or loops.
Let us first introduce some notation. Let $c$ be the distance on $E(G)$
edges of~$G$ defined as follows: for all $e$, $e' \in E(G)$, let
$c(e,e')$ be the minimum number of edges necessary to connect $e$ to
$e'$ plus one ($0$ if
$e=e'$, and $1$ if $e\sim e'$). For all $e=\{i,j\}$, let $\C(2,e)$
be the set of maximal complete $2$-partite subgraphs $S\subset G$ such
that $i, j \in S$ and, for all $k, l \in S$ with $k\sim l$,
$a_{k,l}=a_{i,j}$.
\begin{prop}
\label{prop:str2}
Assume the graph does not contain triangles nor loops. If, for some
$e\in E(G)$,
%
\begin{equation}
\label{eq:2part}
\min_{S\in\C(2,e)}\max_{k\in S, l\in\partial S}a_{k,l}\le a_e,
\end{equation}
then there exists at least one strongly trapping subset.

Note that (\ref{eq:2part}) holds if
\[
\max_{c(e,e')\le2}a_{e'}\le a_e.
\]
\end{prop}
\begin{rem}
\label{rem:str2}
If, for all $e\in E(G)$, (\ref{eq:2part}) does not hold, then there
exists, for all $e\in E(G)$, an infinite sequence of edges $(e_n)_{n\in
\N_0}$ such that $e_0=e$, $e_n\sim e_{n+1}$ and, for all $n\in\N$,
$a_{e_{n}}\le a_{e_{n+1}}$ and $a_{e_{n}}< a_{e_{n+2}}$. However, even
in this case, there can exist a strictly stable equilibrium $x\in\E
_s$ (but no strongly trapping subset).
\end{rem}

\begin{pf*}{Proof of Proposition \ref{prop:str2}}
By assumption, there
exist $e=\{i,j\}$ and a~maximal complete $2$-partite subgraph $S\subset
G$ containing $i$ and $j$, with partitions $V_1$ and $V_2$, and
satisfying conditions (a), (b) and (c)(i) of
Definition~\ref{def:traps}. For all $k\in\partial S$, $k$ is adjacent
to at most one of two partitions, say, $V_1$, since otherwise $G$
would\vadjust{\goodbreak}
contain a triangle; if $k$ were adjacent to all vertices in~$V_1$, then
it would be in $V_2$, since $S$ is assumed maximal. Hence,
(c)(ii) holds as well, and $S$ is a strongly trapping subset.
\end{pf*}

When the graph contains triangles, the property outlined in Remark \ref
{rem:str2}, that is, the existence of an infinite sequence of edges
with increasing labels when there is no strongly trapping subset, does
not hold anymore. The~maxi\-mum of the Lyapounov function on a complete
subgraph with more than two vertices takes a nontrivial form, which can
lead to counterintuitive behavior.

We show, for instance, in Example \ref{ex:2} a case where the
reinforcement matrix~$a$ has a strict global maximum at a certain edge,
but where, however, there is no stable equilibrium at all. We believe
the walk to be transient in this example.

\begin{ex}
\label{ex:1}
Let us show, in the case $a=(\1_{i\sim j})_{i,j\in G}$, that Theorem
\ref{thm:conv} is stronger than Theorem \ref{thm:traps}. Consider a
graph $G$ on six vertices
$A$, $B$, $C$, $D$, $E$ and $F$, with a neighborhood relation $\sim$
defined as follows (see Figure~\ref{fig1}): $A\sim B\sim C\sim D\sim
A$, $C\sim E\sim D$ and $E\sim F$
(recall that the graph~$G$ is symmetric). Let
$x=(x_A,x_B,x_C,x_D,x_E,x_F):=(3/8,3/8,1/8,1/8,0,0)$, then $S(x)=\{
A,B,C,D\}$ and $\partial S(x)=\{E\}$.
Also, $x$ is
an equilibrium of~(\ref{eqdiff}), (P)$_{S(x)}$ is satisfied
with $V_1=\{A,C\}$, $V_2=\{B,D\}$, and $N_E(x)=1/4<H(x)=1/2$, which
implies that $x$ is a strictly stable equilibrium by Theorem~\ref
{partition2}, hence subsequently by Theorem~\ref{thm:conv} that $\Rc
=T(x)$ with positive probability.

Now let us prove by contradiction that $T(x)$ with such $x$ does not
satisfy the assumptions of Theorem \ref{thm:traps} above. Indeed, if
$T(x)=S\cup\partial S$, then $S\subset\{A,B,C,D\}$ since, otherwise,
$F$ would belong to $T(x)$. Now the condition that, for all $i\in
\partial S$, $\exists p\in\{1,\ldots,d\}$ and $j\in S\setminus V_p$
such that $i\not\sim V_p\cup\{j\}$ implies, in particular, that a
vertex in $\partial S$ is not connected to at least two other vertices
in $S$, so that $i\in\partial S$ cannot be $A$, $B$, $C$ or $D$, which
are connected to all other but one vertex in $\{A,B,C,D\}$. Hence, $S=\{
A,B,C,D\}$, but then $i:=E$ is connected to both partitions of $S$, and
does not satisfy the condition mentioned in the last sentence, bringing
a contradiction.
\end{ex}

\begin{ex}
\label{ex:2}
Let us first study the case of a triangle $(G,\sim)$, $G:=\{0,1,2\}$,
$0\sim1\sim2\sim0$, with reinforcement coefficients $a:=a_{0,1}$,
$b:=a_{1,2}$,  $c:= a_{0,2} >0$.

If $a<b+c$, then the equilibrium $x=(x_0,x_1,x_2)=(1/2,1/2,0)$ is not
stable, since $N_2(x)=(b+c)/2>H(x)=a/2$. Hence, if we assume that
%
\begin{equation}
\label{c:triangle}
a<b+c,\qquad  b<a+c,\qquad  c<a+b,
\end{equation}
then a stable equilibrium has to belong to the interior of the simplex
$\De$. A~simple calculation shows that there is only one such equilibrium:
\[
x=(x_0,x_1,x_2):= \biggl(\frac{c(a+b-c)}{\delta},\frac{b(a+c-b)}{\delta
},\frac{a(b+c-a)}{\delta} \biggr),
\]
where
\[
\delta:=(a+b+c)^2-2(a^2+b^2+c^2);
\]
$\delta\,{>}\,0$, which can be shown by adding up inequalities
$(b-a)^2\,{\le}\,c^2$, \mbox{$(c-a)^2\,{\le}\,b^2$} and $(c-b)^2\,{\le}\,a^2$. Then $H(x)=2abc/\delta$.

\begin{figure}

\includegraphics{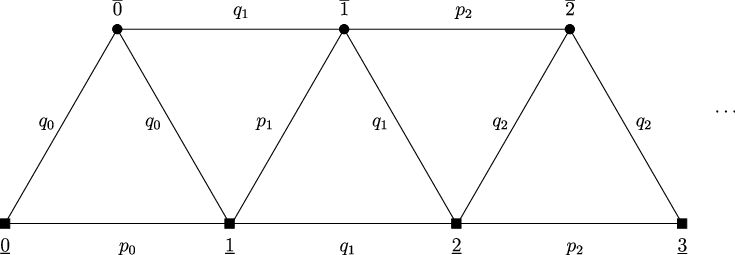}

\caption{On the infinite graph on the figure, with reinforcement
coefficient sequences $(p_n)_{n\ge0}$ strictly decreasing and
$(q_n)_{n\ge0}$ strictly increasing, we show in Example \protect\ref{ex:2}
that, even if $p_0=\sup_{n\ge0}p_n>\sup_{n\ge0}q_n$, we can choose
these sequences in such a way that there is no stable equilibrium in
$\De$, and therefore no trapping subgraph.} \label{fig2}
\end{figure}

Let $\N:=\Z_+$. Let us now consider the following graph $(G,\sim)$
with vertices $G:=\{\ui, \oi, i\in\N\}$ and adjacency $\ui\sim
\uip$, $\oi\sim\oip$, $\ui\sim\oi$ and $\oi\sim\uip$, for all
$i\in\N$, as
drawn in Figure \ref{fig2}.

Fix $\epsilon, \eta, p, q >0$, $\mu\in(0,1)$, which will be
chosen later. Let, for all \mbox{$n\in\N$},
%
\begin{equation}
\label{eq:pnqn}
p_n:=p\prod_{k=0}^{n-1}(1-\mu^k\epsilon),\qquad  q_n:=q\prod
_{k=0}^{n-1}(1+\mu^k\eta).
\end{equation}

Note that, for all $n\in\N$,
\[
p \biggl(1-\frac{\epsilon}{1-\mu} \biggr)\le p_n\le p,\qquad  q\le q_n\le qe^{
{\eta}/{(1-\mu)}}.
\]

Now assume that the reinforcement matrix $(a_{k,l})_{k,l\in G}$ is
defined as follows, depending on $(p_n)_{n\in\N}$ and $(q_n)_{n\in\N
}$, for all $i\in\N$:
\begin{eqnarray*}
a_{\uie,\uio}&:=&p_{2i},\qquad  a_{ \oio, \oiep}:=p_{2(i+1)},\\
a_{\oie,\oio}&\hspace*{3pt}=&a_{\uio,\uiep}:=q_{2i+1}\\
a_{\uie,\oie}&\hspace*{3pt}=&a_{\oie,\uio}:=q_{2i}\\
a_{\uio,\oio}&:=&p_{2i+1}\\
a_{\oio,\uiep}&:=&q_{2i+1}.
\end{eqnarray*}

Let $x\in\De$ be a stable equilibrium of (\ref{eqdiff}). Then, by
Theorem \ref{partition}, (P)$_{S(x)}$ holds, so that $S(x)$\vadjust{\goodbreak}
consists of two vertices or a triangle [it cannot be made of four
vertices, because of (P)$_{S(x)}$(c)]. Assume
%
\begin{equation}
\label{a:1ex2}
p<2q,\qquad  \eta qe^{{\eta}/{(1-\mu)}}<p \biggl( 1-\frac{\epsilon}{1-\mu} \biggr).
\end{equation}
Then, for all $i\in\N$,
\[
p_i<2q_i,\qquad  p_{i+1}<q_i+q_{i+1},\qquad  q_{i+1}<q_i+p_{i+1},
\]
so that $S(x)$ has to be a triangle.

Assume $S(x):=\{\uie,\oie,\uio\}$ for some $i\in\N$; the argument
is similar in other cases. Then

\[
x_{\oio}=\frac{H(x)}{2q_i},\qquad
x_{\oie}=\frac{H(x)}{2q_i^2}(2q_i-p_i),\qquad
x_{\uio}=\frac{H(x)}{2q_i},
\]
and
\[
N_{\uie}(x)=q_ix_{\oie}+p_ix_{\uio}=H(x),
\]
and, therefore,
\begin{eqnarray*}
N_{\oio}(x)&=&q_{i+1}x_{\oie}+p_{i+1}x_{\uio}\\
&=&H(x)+\frac{H(x)}{2q_i^2}[(q_{i+1}-q_i)(2q_i-p_i)+(p_{i+1}-p_i)q_i]\\
&=&H(x)+\frac{H(x)}{2q_i^2}\mu^i[\eta q_i(2q_i-p_i)-\epsilon p_iq_i]>H(x)
\end{eqnarray*}
if
%
\begin{equation}
\label{a:2ex2}
\eta>\epsilon\frac{p}{2q-p},
\end{equation}
using that $p/(2q-p)>p_i/(2q_i-p_i)$ for all $i\in\N$.

Hence, $x$ is not a stable equilibrium, which leads to a contradiction.
\end{ex}
%
\section{Introduction to the proofs}
\subsection{Notation}
\label{sec:notation}
We let $\N:=\Z_+$, $\N^*:=\N\setminus\{0\}$, $\R_+^*:=\R
_+\setminus\{0\}$.

For all $y=(y_i)_{i\in G}\in\R^G$ and for any finite subset $A$ of
$G$, let
\[
y_A:=\sum_{i\in A}y_i.
\]

Given $r\in\N^*$,
let $(\cdot,\cdot)$ (resp., $|\cdot|$, $\|\cdot\|_\iy$) be the
scalar product
(resp.,
the canonical norm, the infinity norm)
on $\R^r$, defined by
\[
(a,b)=\sum_{i=1}^ra_ib_i,\qquad  |a|=\sr{(a,a)},\qquad  \|a\|_\iy:=\max_{1\le
i\le r}|a_i|
\]
if $a=(a_1,\ldots,a_r)$ and\vadjust{\goodbreak} $b=(b_1,\ldots,b_r)$.\vfill\eject

Given a real $r\times r$
matrix $M$ with real eigenvalues, we let $\las(M)$ denote the set of
eigenvalues of $M.$ When $M$ is symmetric we let
$M[\cdot]$ denote the quadratic form associated
to $M$, defined by
$M[a]=(Ma,a)$ for all $a\in\R^r$.

Given $y_1, \ldots, y_r$, we let $\Di(y_1,\ldots,y_r)$ be the diagonal
$r\times r$ matrix of diagonal terms $y_1, \ldots, y_r$.

For all $u, v \in\R$, we write $u=\Box(v)$ if $|u|\le v$. Given
two (random) sequences $(u_n)_{n\ge k}$ and $(v_n)_{n\ge k}$ taking
values in $\R$, we write $u_n\equiv v_n$ if $u_n-v_n$ converges a.s.,
and $u_n\sim_{n\to\iy}v_n$ iff $u_n/v_n\to_{n\to\iy}1$, with the
convention that $0/0=1$.

Let $\Cst(a_1,a_2,\ldots, a_p)$ denote a positive constant depending
only on $a_1, a_2,\allowbreak \ldots, a_p$, and let $\Cst$ denote a universal
positive constant.

\subsection{\texorpdfstring{Proof of Theorems \protect\ref{partition}, \protect\ref{partition2} and \protect\ref{thm:traps}}
{Proof of Theorems 1, 2 and 4}}
\label{sec:spartition}
Theorems \ref{partition} and \ref{partition2} are a consequence of
the more general three following Lemmas \ref{lem:1}, \ref{lem:2} and
\ref{lem:3} below.

\subsubsection{\texorpdfstring{Lemmas \protect\ref{lem:1}, \protect\ref{lem:2} and \protect\ref{lem:3},
and proof of Theorem \protect\ref{partition}}
{Lemmas 1, 2 and 3, and proof of Theorem 1}}

$\!\!\!$By the following~Lem\-ma~\ref{lem:1}, if an equilibrium $x\in\De$ is
stable, then the eigenvalues of $[a_{i,j}-2H(x)]_{i,j\in S(x)}$, which
depend only on $a$, $S(x)$ and $H(x)$, are nonpositive. This property
will subsequently imply (P)$_{S(x)}$, by Lemmas \ref{lem:2} and
\ref{lem:3}.

\begin{lemma}
\label{lem:1}
Let $x=(x_i)_{i\in G}\in\De$ be an equilibrium. Then:
\begin{enumerate}[(a)]
\item[(a)] $DF(x)$ has real eigenvalues.
\item[(b)] The three following assertions are equivalent:
\begin{enumerate}[(iii)]
\item[(i)] $x$ is stable,
\item[(ii)] $\max\las(DF(x))\le0$,
\item[(iii)] $\max(\las( [a_{i,j}-2H(x) ]_{i,j\in S(x)} )\cup\{
N_i(x)-H(x),i\in\partial S(x) \} ) \leq0$.
\end{enumerate}
\item[(c)] If $x$ is stable, then it is feasible.

\end{enumerate}
\end{lemma}

Lemma \ref{lem:2} yields an algebraically simpler
characterization of assertion (P)$_{S}$ for $S\subset G$; recall
that, given subsets $S$ and $R$ of $G$, $\partial_S R$, defined in
Section~\ref{sec:gintro}, is the outer boundary of $R$ inside $S$.

\begin{lemma}
\label{lem:2}
The statement \textup{(P)}$_S$ is equivalent to
\begin{enumerate}[$\mathrm{(P)}_S'$]
\item[$\mathrm{(P)}_S'$]
If $j, k \in S$ are such that $j\not\sim k$, then,
for all $i\in S$, $a_{i,j}=a_{i,k}$
(so that $\partial_{S}\{j\}=\partial_{S}\{k\}$ in particular).
\end{enumerate}
\end{lemma}

Lemma \ref{lem:3} states that (P)$_{S(x)}$ holds if the
eigenvalues of $[a_{i,j}-2H(x)]_{i,j\in S(x)}$ are nonpositive, with
equivalence if $a=(\1_{i\sim j})_{i,j\in G}$.

\begin{lemma}
\label{lem:3}
Let $x=(x_i)_{i\in G}\in\De$ be a feasible equilibrium. Then
\[
\max\las\bigl( [a_{i,j}-2H(x) ]_{i,j\in S(x)} \bigr)\le0\quad
\Longrightarrow\quad \mbox{\textup{(P)}}_{S(x)}'.
\]
If, for some $c>0$, $a_{i,j}=c\1_{i\sim j}$ for all $i, j \in
S(x)$, then the above implication is an equivalence.
\end{lemma}

Lemmas \ref{lem:1}, \ref{lem:2} and \ref{lem:3} are proved,
respectively, in Sections \ref{sec:lem1}, \ref{sec:lem2} and \ref
{sec:lem3}. They obviously imply
Theorem \ref{partition}.

\subsubsection{\texorpdfstring{Proof of Theorem \protect\ref{partition2}}{Proof of Theorem 2}}
\label{p:partition2}
Suppose $a=(\1_{i\sim j})_{i,j\in G}$, and let $x\in\De$.

First assume that $(S(x),\sim)$ contains no loop. If $x$ is a stable
equilibrium, then (P)$_{S(x)}$ and, thus, (i) holds by
Theorem \ref{partition}; let $V_k$, $1\le k\le d$ be the partitions of
$S(x)$. Then $d\ge2$ [otherwise $H(x)=0$ and $x$ is not feasible, thus
not stable by Lemma \ref{lem:1}] and, for
all $1\le k\le d$, $j\in V_k$,
\[
v_k:=\sum_{i\in V_k}x_i= 1-N_j(x)=1-H(x),
\]
so that $v_k=1/d$ (since $\sum_k v_k=1$) and $H(x)=1-1/d$, and,
subsequently, (ii)--(iii) hold by Lemma \ref{lem:1}. Conversely,
assume (i)--(iii) hold; then $N_i(x)=1-1/d$ for all $i\in S(x)$,
so that $H(x)=\sum_{i\in S(x)}x_iN_i(x)=1-1/d$ and $x$ is a~feasible
equilibrium. Now (i) implies (P)$_{S(x)}$ and thus
(P)$_{S(x)}'$ by Lemma~\ref{lem:2}. Hence, using Lemmas \ref{lem:1}
and \ref{lem:3}, $x$ is a stable equilibrium.

Now assume on the contrary that $(S(x),\sim)$ contains one loop $i\sim
i$. If $x$ is a stable equilibrium, then (P)$_{S(x)}$ again holds
by Theorem \ref{partition}: (P)$_{S(x)}$(b) implies
$N_i(x)=1=H(x)$ ($x$ equilibrium). Hence,
for all $j\in S(x)$, $N_j(x)=1$ and $j\sim k$ for all $k\in S(x)$, that
is, $(S(x),\sim)$ is a clique of loops. Conversely, if $(S(x),\sim)$
is a clique of loops, then (P)$_{S(x)}$ obviously holds so that,
by Lemmas \ref{lem:1} and \ref{lem:3}, $x$ is stable [since $H(x)=1$,
then $N_i(x)\le H(x)$ for all~\mbox{$i\in G$}].

\subsubsection{\texorpdfstring{Proof of Theorem \protect\ref{thm:traps}}{Proof of Theorem 4}}
\label{sec:thm-traps}
First observe that
\[
\Sigma=\Sc(S)\cap\Ec_s.
\]
Indeed, the proof of Theorem \ref{partition2} implies that $\Sigma
\supset\Sc(S)\cap\Ec_s$ and, conversely, that if $x\in\Sigma$,
then $x$ is a equilibrium and, by (c)(ii), for all $j\in\partial
S(x)$, $N_j(x)<H(x)$
[$=a_S(1-1/d)$ if $(S(x),\sim)$ contains no loops, $=a_S$ otherwise],
using assumptions (a)--(b) and (c)(ii) or the second part of
(c)$'$. Also, (P)$_{S(x)}$ holds by (c) or (c)$'$, and,
therefore, $x$ is strictly stable by Lemmas \ref{lem:1}--\ref{lem:3}.
The rest of the proof follows from Theorem \ref{thm:conv}.

\subsection{\texorpdfstring{Proof of Theorem \protect\ref{thm:conv}}{Proof of Theorem 3}}
\label{sec:thm-conv}
First, we provide in Lemma \ref{bound_rn} a rigorous mathematical
setting for the stochastic approximation of the density of occupation
of the VRRW $x(n)$ by solutions of the ordinary differential equation
(\ref{eqdiff}) on a~finite graph~$G$, heuristically justified in
Section \ref{sec:gintro} [see (\ref{mdyn})]. Second, we make use of
this technique and of an entropy function originally introduced in
\cite{losert} to study the VRRW on the finite subgraph $T(x)$ when its
density of occupation is in the neighborhood of a strictly stable
equilibrium $x$, in Lemmas~\ref{lem:evolent}--\ref{lem:aborder}.
Third, we focus again on a general graph $G$---possibly infinite---and
prove in Proposition \ref{prop:conv}, assuming again that the density
of occupation is in the neighborhood of an element $x\in\E_s$, that
the walk eventually localizes in~$T(x)$ with lower bounded\vadjust{\goodbreak} probability.

In the first step, we make use of a technique originally introduced by
M\'{e}tivier and Priouret in 1987 \cite{metivier} and adapted by
Bena{\"{\i}}m \cite{benaim1} in the context of vertex reinforcement
when the graph is complete (Hypothesis 3.1 in \cite{benaim1}).
In Sections~\mbox{\ref{stoc-app}--\ref{sec:bound rn}}, we generalize it and
show that a certain quantity $z(n)$, depending only on $a$, $x(n)$,
$X_n$ and $n$ and defined in (\ref{defzn}), satisfies the recursion
(\ref{eq stoc-approx}):
\[
z(n+1)=z(n)+\frac{1}{n+n_0+1}\frac{F(z(n))}{H(x(n))}+\epsilon_{n+1}+r_{n+1},
\]
where $\Es(\epsilon_{n+1} | \F_n)=0$. The following Lemma \ref
{bound_rn}, proved in Section \ref{sec:bound rn}, provides upper
bounds on the infinity norms of $\epsilon_{n+1}$, $r_{n+1}$ and
$z(n)-x(n)$, and on the conditional variances of $(\epsilon_{n+1})_i$,
$i\in G$.

More precisely, let us break down the set of vertices of $G$ as
$G=S\cup\partial S$, where $(S,\sim)$ is finite, connected and not a
singleton unless it is a loop. Let, for all $\al\in\R_+\setminus\{
0\}$,
%
\begin{equation}
\label{eq:lalph}
\La_\al:=\{x=(x_j)_{j\in G}\in\De\mbox{ s.t. } x_j\ge\al\mbox{
for all } j\in S\}.
\end{equation}

\begin{lemma}
\label{bound_rn}
For all $n\ge\Cst(\al)$ and $i\in G$, if $x(n)\in\La_\al$, then
\begin{eqnarray*}
\mbox{\textup{(a)}}\quad\!\!\|\epsilon_{n+1}\|_\iy&\le&\frac{\Cst
(\al,a,|G|)}{n+n_0},\qquad
\mbox{\textup{(b)}}\quad\!\!\Es((\epsilon_{n+1})_i^2 | \F_n)\le
\frac
{\Cst(\al,a,|G|)x(n)_i}{(n+n_0)^2},\\
\mbox{\textup{(c)}}\quad\!\!\|r_{n+1}\|_\iy&\le&\frac{\Cst(\al
,a,|G|)}{(n+n_0)^2},\qquad
\mbox{\textup{(d)}}\quad\!\!\|z(n)-x(n)\|_\iy\le\frac{\Cst(\al
,a,|G|)}{n+n_0}.
\end{eqnarray*}
\end{lemma}

Note that if $G$ were a complete $d$-partite finite graph for some
$d\ge1$ or, more generally, if $G$ were without loop and, for all
$i, j \in G$ with $i\sim j$, \mbox{$\{i,j\}\cup\partial\{i,j\}=G$}, then
the constants in the inequalities of Lemma \ref{bound_rn} would not
depend on $\al>0$ and, as a consequence, the stochastic approximation
of~$z(n)$ by (\ref{eqdiff}) would hold uniformly a.s.
Indeed, for all $n\in\N$, by the pigeonhole principle, there exists
at least one edge $\{i,j\}$ $i, j\in G$, $i\sim j$, on which the
walk has spent more than $n/|G|^2$ times, so that $x(n)_i\wedge
x(n)_j\ge\frac{1}{|G|^2}\frac{n}{n+n_0}$ and, under the assumption
on $G$, Lemma \ref{bound_rn} with $S:=\{i,j\}$ would yield the claim.

In the second step, we define an entropy function $V_q(\cdot)$,
measuring a
``distance'' between $q$ and an arbitrary point [as can be seen by
(\ref
{eqVqiy2}) below], originally introduced by Losert and Akin in 1983 in
\cite{losert} in the study of the deterministic Fisher--Wright--Haldane
population genetics model, and to our knowledge so far only used for
the analysis of deterministic replicator dynamics. Note that it is not
mathematically a distance, however, since it does not satisfy the
triangle inequality in general.

In the following, until after the statement Lemma \ref{lem:aborder}---and, in particular,
in Lemmas \ref{lem:evolent}--\ref{lem:aborder}---we assume that $x\in\Ec_s$ and $G=T(x)=S(x)\cup\partial S(x)$; this
choice will be justified later in the proof. Note that if $q\in\Nn
(x)\cap\Ec_s$, whe\-re~$\Nn(x)$ is an adequately chosen neighborhood
of $x$, then $q\in\Sc(S(x))$ since $x\in\Ec_s$, so that
$T(q)=T(x)$. Set $S:=S(x)$, $T:=T(x)$, and $\Sc:=\Sc(S(x))$ for simplicity.

Lemmas \ref{lem:evolent} and \ref{lem:hestimates} below will imply
that, given any stable equilibrium $q\in\Nn(x)\cap\Ec_s$ as a
reference point, $V_q(z(n))$ decreases in average when $z(n)$ is close
enough to~$x$. Therefore, martingale estimates will enable us to prove
in Lemma \ref{lem:posprob} that, starting in the neighborhood of $x$,
$x(n)$ remains close to~$x$ with large probability if $n$ is large, and
converges to one of the strictly stable equilibria in this neighborhood.

For all $q=(q_i)_{i\in G}\in\Sc$ and $y\in\R^{G}$, let
\[
V_q(y):=
\cases{
\displaystyle -\sum_{i\in S}q_i\log(y_i/q_i )+2y_{\partial S}, &\quad  if $y_i>0$,
$\forall i\in S$,\cr
\iy,&\quad  otherwise.
}
\]

Let, for all $q\in\Sc$ and $r>0$,
\[
B_{V_q}(r):=\{y\in\De\mbox{ s.t. } V_q(y)<r\},\qquad  B_\iy(q,r):=\{y\in
\De\mbox{ s.t. }\|y-q\|_\iy<r\}.
\]

Then, we will prove in Section \ref{sec:evolent} that, for all $q\in
\Sc$, there exist increasing continuous functions $u_{1,q}$,
$u_{2,q}\dvtx\R_+\longrightarrow\R_+$ such that
$u_{1,q}(0)=u_{2,q}(0)=0$ and,
for all $r>0$,
%
\begin{equation}
\label{eqVqiy2}
B_\iy(q,u_{1,q}(r))\subset B_{V_q}(r)\subset B_\iy(q,u_{2,q}(r)).
\end{equation}

Let, for all $q, z \in\R^{G}$,
%
\begin{equation}
\label{eq:iq}
I_q(z):=-\sum_{i\in S}q_i[N_i(z)-H(z)]+2\sum_{i\in\partial
S}z_i[N_i(z)-H(z)].
\end{equation}

The following Lemma \ref{lem:evolent}, also proved in Section \ref
{sec:evolent}, provides the stochastic approximation equation for
$V_q(z(n))$, $q\in\Sc\cap\Ec_s$; we make use of notation $u=\Box
(v)\iff|u|\le v$, introduced in Section \ref{sec:notation}.

\begin{lemma}
\label{lem:evolent}
Let $q\in\Sc\cap\Ec_s$. There exist an adapted process $(\ze
_n)_{n\in\N}$ (not depending on $q$ and $a$), and constants $n_1$ and
$\epsilon$ (depending only on $q$ and $a$) such that, if $n\ge n_1$
and $x(n)\in B_{V_q}(\epsilon)$, then $V_q(z(n))<\iy$, \mbox{$V_q(z(n+1)) <\iy$}, $\Es(\zeta_{n+1} | \F_n)=0$ and
\begin{eqnarray}
\label{evolent}
V_q\bigl(z(n+1)\bigr) &=&
V_q(z(n))+\frac{I_q(x(n))}{(n+n_0+1)H(x(n))}\nonumber\\[-8pt]\\[-8pt]
&&{} -(q,\ze_{n+1})+2(\epsilon_{n+1})_{\partial S}+\Box\biggl(\frac{\Cst
(q,a)}{(n+n_0)^2} \biggr).\nonumber
\end{eqnarray}
\end{lemma}

Lemma \ref{lem:hestimates}, proved in Section \ref
{sec:hestimates}, provides estimates of the Lyapounov function~$H$, and
of $I_\cdot(\cdot)$, in the neighborhood of a strictly stable equilibrium. It
will not only be useful in the proof of Lemma \ref{lem:posprob},
stating convergence of $x(n)$ with large probability, but also for
Lemma \ref{lem:speed} on the rate of this\vadjust{\goodbreak} convergence.
\begin{lemma}
\label{lem:hestimates}
There exists a neighborhood $\Nn(x)$ of $x$ in $\De$ such that, for
all $q\in\Nn(x)\cap\Ec_s$, $y\in\Nn(x)$,
%
\begin{eqnarray}
\label{in:lem4}
&&\mbox{\textup{(a)}}\quad \Cst(x,a)J(y)\le H(q)-H(y)\le\Cst
(x,a)J(y),
\\
\label{Iq<0}
&&\mbox{\textup{(b)}}\quad {-}[H(q)-H(y)+\Cst(x,a)y_{\partial
S}]\nonumber\\[-8pt]\\[-8pt]
&&\hphantom{\mbox{\textup{(b)}}\quad}\qquad \le
I_q(y)\le-[H(q)-H(y)+\Cst(x,a)y_{\partial S}]
\le0.\nonumber
\end{eqnarray}
\end{lemma}

\begin{rem}
Lemma \ref{lem:hestimates} implies that $y\in\Nn(x)$ is an
equilibrium iff $H(y)=H(x)$. Also note that the maximality of $H$ at
$x\in\Ec_s$ is not global in general. For instance, in the
counterexample at the end of Section \ref{sec:gintro},
$x:=(3/8,3/8,1/8,1/8, 0)\in\Ec_s$, but, letting
$y:=(0,0,1/3,1/3,1/3)$, $H(y)=2/3>H(x)=1/2$.
\end{rem}

The proof of Lemma \ref{lem:posprob} is shown in Section \ref
{sec:enttrap}. A key point in its proof is that the martingale term
$-(q,\ze_{n+1})+2(\epsilon_{n+1})_{\partial S}$, in Lemma \ref
{lem:evolent}, is a linear function of $\ze_{n+1}$ and $\epsilon
_{n+1}$ which do not depend on $q$, so that the two corresponding
convergence results of these martingales will apply from any reference point
$q\in\Ec_s\cap\Nn(x)$. It will enable us to prove that, if $r$ is a
accumulation point of $x(n)$, then $V_r(x(n))$ a.s. converges to $0$ if
$r\in\Nn(x)$ although~$r$ is random.

\begin{lemma}
\label{lem:posprob}
There exist $\epsilon_0:=\Cst(x,a)$ and $n_1:=\Cst(x,a)$ such that,
if for some $\epsilon\le\epsilon_0$ and $n\ge n_1$, $x(n)\in
B_{V_x}(\epsilon/2)$, then
\[
\Pb\bigl(\Lc(B_{V_x}(\epsilon)) | \F_n\bigr)\ge1-\exp\bigl(-\epsilon^2\Cst
(x,a)(n+n_0)\bigr).
\]
\end{lemma}

Next, we provide in Lemma \ref{lem:speed} some
information on the rate of convergence of $x(n)$ to $x(\iy)$, which
will be necessary for the asymptotic estimates on the frontier $\A
_\partial(x(\iy))$ in Lemma \ref{lem:aborder}.

\begin{lemma}
\label{lem:speed}
There exist $\epsilon$, $\nu :=\Cst(x,a)$ such that, a.s. on $\Lc
(B_{V_x}(\epsilon))$,
\[
\lim_{n\to\iy}\bigl(x(n)-x(\iy)\bigr)n^{\nu}=0.
\]
\end{lemma}

The proof of Lemma \ref{lem:speed}, given in Section \ref{sec:speed},
starts with a preliminary estimate of the rate of convergence of
$H(x(n))$ to $H(x(\iy))$. To this end, we make use of Lemma \ref
{lem:stoc h} below, giving the stochastic approximation equation of
$H(z(n))$. It implies, together with Lemma \ref{lem:hestimates}(a), that the
expected value of $H(z(n+1))-H(z(n))$ is at least $\Cst
(x,a)(H(x)-H(z(n))$, so that we can then estimate the rate of $H(x(n))$
to $H(x)$ by a one-dimensional technique.

Finally, this estimate implies similar ones for the convergence of
$J(x(n))$ and $I_{x(\iy)}(x(n))$ to $0$ by Lemma \ref
{lem:hestimates}, so that we conclude using entropy estimates for the\vadjust{\goodbreak}
rate of convergence of $V_{x(\iy)}(z(n))$, using again that only two
martingales estimates are necessary, given the linearity of the
perturbation in (\ref{evolent}) with respect to the reference point
$q\in\Ec_s\cap\Nn(x)$.

\begin{lemma}
\label{lem:stoc h}
For all $n\in\N$,
%
\begin{equation}
\label{stoc_h}
H\bigl(z(n+1)\bigr)-H(z(n))=\frac{1}{n+n_0+1}\frac{J(z(n))}{H(x(n))}+\xi_{n+1}+s_{n+1},
\end{equation}
where $\Es(\xi_{n+1} | \F_n)=0$ and, if for some $\al>0$, $x(n)\in
\La_\al$ and $n\ge\Cst(\al)$, then
\[
(1)\quad\|\xi_{n+1}\|_\iy\le\frac{\Cst(\al
,a,|G|)}{n+n_0},\qquad
(2)\quad\|s_{n+1}\|_\iy\le\frac{\Cst(\al,a,|G|)}{(n+n_0)^2}.
\]
\end{lemma}
Lemma \ref{lem:stoc h} is proved in Section \ref{sec:stoc h}.

Lemma \ref{lem:aborder} yields the asymptotic behavior on the
border sites $\partial S$. This behavior is similar to the one one
would obtain without perturbation [i.e., with $(\epsilon_{n})_{n\in\N
^*}=0$ in (\ref{eq stoc-approx})]. Indeed, if $i\in\partial S$, then
$N_i(x)-H(x)<0$ is the eigenvalue of the Jacobian matrix of (\ref
{eqdiff}) in the direction $(\delta_{i,j})_{j\in G}$
(see the proof of Lemma~\ref{lem:1}), and the renormalization in time
is approximately in $H(x)^{-1}\log n$ [see equation~(\ref{eq
stoc-approx})], so that the replicator equation (\ref{eqdiff}) would
predict that $i\in\partial S$ is visited of the order of
$n^{N_i(x)/H(x)}$ times at time $n$. This similarity with the noiseless
case is due to the fact that the perturbation $(\epsilon_{n})_{n\in\N
^*}$ is weak near the boundary [see Lemma~\ref{bound_rn}(b)].

\begin{lemma}
\label{lem:aborder}
There exists $\epsilon:=\Cst(x,a)$ such that, a.s. on $\Lc
(B_{V_x}(\epsilon))$, $\A_\partial(x(\iy))$ occurs a.s.
\end{lemma}

The proof of Lemma \ref{lem:aborder}, given in Section \ref
{sec:aborder}, makes use of a martingale technique developed in \cite
{tarres}, Section 3.1, and in \cite{limtar} in the context of strong
edge reinforcement. We could have shown Lemma \ref{lem:aborder}
by a thorough study of the border sites coordinates of the stochastic
approximation equation (\ref{eq stoc-approx}), but it would lead to a
significantly longer---and less intuitive---proof.

Now we do not assume anymore that $G=T(x)$ for some $x\in\De$, in
other words, we let the graph $(G,\sim)$ be arbitrary, possibly infinite.

Let, for all $n, k \in\N\cup\{\iy\}$, $n\ge k$, $\Rc_{n,k}$ be
the range of the vertex-reinforced random walk between times $n$ and
$k$, that is,
\[
\Rc_{n,k}:=\{i\in G\mbox{ s.t. } X_j=i\mbox{ for some }j\in[n,k]\};
\]
note that, for all $n\in\N$, $\Rc\subset\Rc_{n,\iy}$.

\begin{prop}
\label{prop:conv}
Let $x\in\Ec_s$. There exists $\epsilon:=\Cst(x,a)$ such that, for
all $n\ge\Cst(x,a)$, if $X_n\in T(x)$ and $x(n)\in B_{V_x}(\epsilon/2)$, then
\[
\Pb\bigl(\{\Rc_{n,\iy}=T(x)\}\cap\Lc(B_{V_x}(\epsilon))\cap\A
_\partial(y(\iy)) | \F_n\bigr)>0.
\]
Moreover, the rate of convergence is at least reciprocally polynomial,
that is, there exists $\nu :=\Cst(x,a)$ such that, a.s. on $\Lc
(B_{V_x}(\epsilon))$,
\[
\lim_{k\to\iy}\bigl(y(k)-y(\iy)\bigr)k^{\nu}=0.
\]
\end{prop}

Proposition \ref{prop:conv} is proved in Section \ref{sec:loc}. It
obviously implies Theorem \ref{thm:conv}: indeed, given  $U$ a
neighborhood of $x$, there exists $\epsilon>0$ such that
$B_{V_x}(\epsilon)\subset U$, and $X_n\in T(x)$ and $x(n)\in
B_{V_x}(\epsilon/2)$ occurs with positive probability if $n$ is large enough.

Observe that, if $G=T(x)$, then this Proposition \ref{prop:conv} is a
direct consequence of Lemmas \ref{lem:posprob}, \ref{lem:speed} and
\ref{lem:aborder}. The localization with positive probability in this
subgraph $T(x)$ results from a Borel--Cantelli type argument: the
probability to visit $\partial T(x)$ at time $n$ starting from
$S(x)$ is, by Lemma \ref{lem:aborder}, upper bounded by a term smaller
than $n^{\al-2}$, where $\al\approx\max_{i\in\partial
S}N_i(x)/H(x)<1$, and $\sum_{n\in\N} n^{\al-2}<\iy$. Technically,
the proof is based on a comparison of the probability of arbitrary
paths remaining in $T(x)$ for the VRRWs defined, respectively, on the
graphs $T(x)$ and $G$.

\subsection{Contents}
Section \ref{sec:detdyn} concerns the results on the deterministic
replicator dynamics: Lemmas \ref{lem:1}--\ref{lem:3} and Lemma \ref
{lem:hestimates} are proved, respectively, in Sections \ref{sec:lem1}--\ref{sec:lem3} and \ref{sec:hestimates}.

Section \ref{sec:stoc} develops the framework relating the behavior of
the vector of density of occupation $x(n)$ to the replicator equation
(\ref{eqdiff}): we write the stochastic approximation equation (\ref
{eq stoc-approx}) in Section \ref{stoc-app}, establish in Section~\ref
{as-trap} some preliminary estimates on the underlying Markov Chain
$M(x)$, prove Lemma~\ref{bound_rn} in Section~\ref{sec:bound rn},
prove Lemmas~\ref{lem:evolent} and~\ref{lem:stoc h} [stochastic
approximation equations for $V_q(z(n))$ and $H(z(n))$] and inclusions~(\ref{eqVqiy2})
in Sections~\ref{sec:evolent} and~\ref{sec:stoc h}.

Section \ref{sec:asymp} is devoted to the proofs of the asymptotic
results for the VRRW: Lemma \ref{lem:posprob} in Section \ref
{sec:enttrap} on the convergence of $x(n)$ with positive probability,
Lemma \ref{lem:speed} in Section \ref{sec:speed} on the corresponding
speed of convergence, Lemma~\ref{lem:aborder} in Section \ref
{sec:aborder} on the asymptotic behavior of the number of visits on the
frontier of the trapping subset, and Proposition \ref{prop:conv} in
Section \ref{sec:loc} on localization with positive probability in the
trapping subsets.

Finally, we show in Appendix \ref{sec:rsm} a lemma on the remainder of
square-bounded martingales, which is useful in the proofs of Lemma \ref
{lem:speed} and Proposition \ref{prop:cex}, whereas Appendix \ref
{sec:cex} is devoted to the proof of Proposition \ref{prop:cex}.

\section{Results on the replicator dynamics}\label{sec:detdyn}

\subsection{\texorpdfstring{Proof of Lemma \protect\ref{lem:1}}{Proof of Lemma 1}}\label{sec:lem1}

Note that $DF(x)v=-H(x)v=0$ if $S(v)\cap T(x)=\varnothing$, so that it is
sufficient to study the eigenvalues of $DF(x)$ on
$\{v\in\R^{G}\mbox{ s.t. } S(v)\subset T(x)\}$; hence, we can assume
that $G$ is finite [equal to $T(x)$] w.l.o.g.\vadjust{\goodbreak}

Let $S:=S(x)$ for convenience.
For all $i, j \in G$,
\[
\frac{\partial F_i}{\partial x_j}
=
\cases{
N_i(x)-H(x),&\quad if $x_i=0$ and $j=i$,\cr
0,&\quad if $x_i=0$ and $j\not=i$,\cr
x_i[a_{i,j}-2H(x)],&\quad if $ x_i\not=0$ and $ x_j\not=0$,\cr
x_i[a_{i,j}-2N_j(x)],&\quad if $ x_i\not=0$ and $ x_j=0$.}
\]
Let us now consider matrix $DF(x)$
by taking the following
order on the indices: we take first the indices $i, j \in
G\setminus S$, and second the indices
$i, j\in S$,
\[
\pmatrix{
\Di\bigl(N_i(x)-H(x)\bigr)_{i\in G\setminus S}& (0)\cr
(*)& DB
},
\]
where
\[
B= [a_{i,j}-2H(x) ]_{i,j\in S},\qquad  D=\Di(x_i)_{i\in S}.
\]

The matrix $DB$ is easily seen to be self-adjoint with respect to the
scalar product $(u,v)_{D^{-1}}:=(D^{-1}u,v).$
Hence, $DB$ has real eigenvalues. This proves the first statement of
the lemma.

Note that if we consider (\ref{eqdiff}) as a differential equation
on $\R^{G}$, then
\[
(F(x),\1)=\frac{d(x(t), \1)}{dt}\bigg|_{t=0,x(0)=x}=-\bigl((x,\1)-1\bigr)H(x).
\]
Therefore, if $x\in\De$ [which implies $(x,\1)=1$], for all vector
$u\in\R^{G}$,
%
\begin{equation}
\label{DF1e}
(DF(x) u , \1) =-H(x) (u, \1).
\end{equation}
Hence, $p\dvtx u\mapsto(u,\1)$ is an eigenvector of $^t DF(x)$ with
eigenvalue $-H(x).$ This makes $-H(x)$ an eigenvalue of $ DF(x)$ and,
more precisely,
\[
\las( DF(x))=\{-H(x)\}\cup\las( DF(x)|_{T\De});
\]
indeed, by (\ref{DF1e}), an eigenvector $u$ of $DF(x)$ with eigenvalue
$\lambda\not=-H(x)$ belongs to $\mathsf{Ker } p=T\De$.
Therefore, the stability of an equilibrium $x$ of (\ref{eqdiff}) on~$\R^{G}$
is equivalent to the stability restricted on $\De$, which
completes the proof of the first equivalence in statement (b).

\begin{claim*}
Let $M=\Di(y_1,\ldots,y_r)$ be a diagonal $r\times r$ matrix,
with $y_1,\ldots,\allowbreak y_r \in\R_+^*$, and
let $N$ be a symmetric $r\times r$ matrix.
Then $\min\las(N)\ge0\iff\min\las(MN)\ge0$ and, under this assumption,
\[
\min\las(MN)\ge\min\las(N)\min\{y_i\}_{1\le i\le r}.
\]
\end{claim*}
\begin{pf} 
It suffices to prove that $\min\las(N)\ge0$ implies $\min\las
(MN)\ge0$ and the
corresponding inequality, since the coinverse
statement is symmetrical.

Recall that, for any $r\times r$ symmetric matrix $R$ with nonnegative
eigenvalues, there exist a diagonal matrix $D$ and an orthogonal matrix
$Q$ such that $R=Q^TDQ$, hence,
\[
\min\las(R)=\inf_{|t|\ge1}(Dt,t)=\inf_{|t|\ge1}(DQt,Qt)=\inf
_{|t|\ge1}(Rt,t).
\]

Let us define $L=\Di(\sr{y_1},\ldots,\sr{y_r})$. Observe that
$L^2=M$. Now $MN=L(LNL)L^{-1}$ implies $\las(MN)=\las(LNL)$.

$LNL$ is symmetric; therefore,
\begin{eqnarray*}
\min\las(MN)&=&\min\las(LNL)=\inf_{|t|\ge1}(LNLt,t)\\
&=&
\inf_{|t|\ge1}(NLt,Lt)
\ge\inf_{|u|\ge\min_{1\le i\le r}\sr{y_i}}(Nu,u)\\
&=&\min_{1\le i\le r} y_i \inf_{|u|\ge1}(Nu,u)=
\min_{1\le i\le r} y_i\las(N).
\end{eqnarray*}
\upqed
\end{pf}

To complete the proof of statement (b), we apply the claim to $M:=D$ and $N:=-B$.

It remains to prove that a stable equilibrium in $\De$ is feasible.
Let \mbox{$x \in\Delta$} be such an equilibrium. Assume that $H(x) = 0$. If $x_i = 0$ for some $i$ then,
by Lemma \ref{lem:1}(b),
$N_i(x) = 0$, so that $x_j = 0$ for all $j \sim i$. Hence, $x = 0$, which is
contradictory. Now, if $x_i \neq0$ for all $i$, then $G$ is necessarily
finite (by definition of $\De$), and
$a=(a_{i,j})_{i,j\in G}= 0$ since its eigenvalues are nonpositive
[Lemma \ref{lem:1}(b) again] and its trace is nonnegative.
This is again contradictory.\looseness=1

\subsection{\texorpdfstring{Proof of Lemma \protect\ref{lem:2}}{Proof of Lemma 2}}
\label{sec:lem2}
Let $\partial:=\partial_S$, $\mathrm{(P)}:=\mathrm{(P)}_S$ and $
\mathrm{(P)}':=\mathrm{(P)}_S'$ for simplicity.

Assume (P) holds for some $d\ge1$.
Let us prove that, if $i, j, k \in S$ are such that $i\sim j\not
\sim k$, then $a_{i,j}=a_{i,k}$.

If $i=j$, then $i=j\not\sim k$ implies, by (P)(a)--(b), that
$k\notin S$---and therefore a~contradiction---since if $k$ were in
$S$, it would be in the partition of $i$, which is a singleton.
If $i\not=j\not\sim k$, then $j$ and $k$ are in the same partition of
$S$. Hence, $a_{i,j}=a_{i,k}$ by (P)(c), which completes the
proof of (P)$'$.

Assume now (P)$'$. Let us prove that the relation $R$ defined on
$S$ by
\[
iRj\quad \Longleftrightarrow\quad  i\not\sim j\mbox{ or } i=j
\]
is an equivalence relation on $S$.
It is clearly
symmetric and reflexive. Let us prove that it is
transitive: let $i, j, k \in S$ be such that $iRj$ and $jRk$, and
prove $iRk$.
This is immediate if $i=j$ or $j=k$; hence, assume that $i\not=j$ and
$j\not=k$; then (P)$'$ implies $\partial_S\{i\}=\partial_S\{j\}
=\partial_S\{k\}$. If we had
$i\sim k$, then it would imply $k\in\partial_S\{i\}=\partial_S\{j\}
$, and, therefore,
$j\sim k$, which leads to a~contradiction.

Now let us prove that there is only one element in the partition of a
loop. Assume that $iRj$, $i\sim i$ and $j\not=i$ for $i, j \in S$;
(P)$'$ implies in this case that
$a_{i,i}=a_{i,j}>0$, so that $i\sim j$, hence, $i=j$ since $iRj$ holds,
which leads to a~contradiction.

Let $V_p$, $p=1, \ldots,d$ be the partitions of $R$: elements of
different partitions are connected, by definition, and (P)(a)--(b)
holds for some $d\ge1$. Let us prove (P)(c): let $p, q \in\{
1,\ldots,d\}$ be such that $p\not=q$, and assume $i\in V_p$, $j\in
V_q$. Let
\[
W_{i,j}:=\{(i',j')\in S^2\mbox{ s.t. }a_{i',j'}=a_{i,j}\}.
\]
By applying (P)$'$ twice, we first obtain that $W_{i,j}\supseteq\{
i\}\times V_q$, and second that $W_{i,j}\supseteq V_p\times V_q$, which
enables us to conclude.

\subsection{\texorpdfstring{Proof of Lemma \protect\ref{lem:3}}{Proof of Lemma 3}}
\label{sec:lem3}
Let $S:=S(x)$ and $\mathrm{(P)}':=\mathrm{(P)}_{S(x)}'$ for
simplicity. Let
\[
B=[a_{i,j}-2H(x)]_{i,j\in S}.
\]
Now  $\max\las(B)\le0\iff \forall t\in\R^S, B[t]\le0$. Observe
that, for all $t=\break(t_i)_{i\in S}\in\R^S$,
\[
B[t]=\sum_{i,j\in S}\bigl(a_{i,j}-2H(x)\bigr)t_it_j=H(t)-2H(x) \biggl(\sum_{i\in S}
t_i \biggr)^2.
\]
Let us assume that (P)$'$ does not hold, and deduce that $B[t]>0$
for some $t\in\R^S$, which will prove the first statement.

There exist $i, j, k \in S$ such that $j\not\sim k$ and
$a_{i,j}\not=a_{i,k}$
[otherwise (P)$'$ would be satisfied]. Let, for all $\lambda\in\R$,
\[
t_\lambda:=\bigl(\1_{\{v=i\}}+\lambda\1_{\{v=j\}}-(1+\lambda)\1_{\{v=k\}
}\bigr)_{v\in S}\in\R^S,
\]
then
\[
B[t_\lambda]\ge2\lambda(a_{i,j}-a_{i,k})-2a_{i,k},
\]
so that $B[t_\lambda]>0$ for some $\lambda\in\R$, which yields the
contradiction.

Let us now assume that (P)$'$ holds, and that $a_{i,j}=c\1_{i\sim
j}$, with $c=1$ for simplicity. First assume $S$ contains no loop.
Then, by Lemma \ref{lem:2}, $S$ is a $d$-partite subgraph for some
$d\ge1$ [$\mathrm{(P)}_S\mathrm{(a)}$ holds]; let $V_1, \ldots,
V_d$ be
its partitions, then
\begin{eqnarray*}
B[t]&=&\sum_{i,j\in S}\bigl(\1_{i\sim j}-2H(x)\bigr)t_it_j
=-2H(x) \biggl(\sum_{i\in S} t_i \biggr)^2+\sum_{i,j\in S}\1_{i\sim j}
t_it_j\\
&=&-2H(x) \biggl(\sum_{k=1}^d v_k \biggr)^2+ \biggl(\sum_{k=1}^dv_k \biggr)^2-\sum_{k=1}^dv_k^2,
\end{eqnarray*}
where, for all $i\in\{1,\ldots,d\}$, $v_k=\sum_{i\in V_k}t_i$.
Therefore,
\[
B[t]=-\bigl(2H(x)-1\bigr) \biggl(\sum_{k=1}^dv_k \biggr)^2
-\sum_{k=1}^dv_k^2\le0,
\]
where we use the fact that $H(x)\ge1/2$, since $H(x)=1-1/d$ and $d\ge
2$ (see proof of Theorem \ref{partition2}, Section \ref{p:partition2}).

Now assume that $S$ contains one loop; then, again by the proof of
Theorem~\ref{partition2}, Section~\ref{p:partition2}, it is a clique
of loops and $H(x)=1$; thus,
\[
B[t]=-2 \biggl(\sum_{i\in S} t_i \biggr)^2+ \biggl(\sum_{i\in S} t_i \biggr)^2=- \biggl(\sum_{i\in
S} t_i \biggr)^2\le0.
\]

\subsection{\texorpdfstring{Proof of Lemma \protect\ref{lem:hestimates}}{Proof of Lemma 6}}
\label{sec:hestimates}
Let us first prove (a) in the case $q:=x$, which will imply
$H(q)=H(x)$ for any equilibrium $q\in\Nn(x)$ and therefore imply~(a)
in the general case. Let $x\in\Ec_s$, and let $y\in T\De$ be
such that $x+y\in\De$. Let $S:=S(x)$ for simplicity.

Recall that $G=S\cup\partial S$. We have
%
\begin{eqnarray}
\label{eq:h0}
\qquad H(x+y)&=&\sum_{i,j\in G}a_{i,j}(x_i+y_i)(x_j+y_j)
=H(x)+2\sum_{i\in G} N_i(x)y_i+H(y)
\\
&=&H(x)+2\sum_{i\in G} \bigl(N_i(x)-H(x)\bigr)y_i+\sum_{i,j\in
G}\bigl(a_{i,j}-2H(x)\bigr)y_iy_j\nonumber \\
\label{eq:h}
&=&H(x)+2\sum_{i\in\partial S} \bigl(N_i(x)-H(x)\bigr)y_i\nonumber\\[-8pt]\\[-8pt]
&&{}+\sum_{i,j\in
S}\bigl(a_{i,j}-2H(x)\bigr)y_iy_j +\sum_{i\in\partial S}w_i(y)\nonumber \\
\nonumber
&\le& H(x)+2\sum_{i\in\partial S} \bigl(N_i(x)-H(x)\bigr)y_i+\sum_{i\in
\partial S}w_i(y).
\end{eqnarray}
In the third equality, we make use of the identity $\sum_{i\in
G}y_i=0$, whereas in the fourth equality we notice that $N_i(x)=H(x)$
for all $i\in S$ and that the reinforcement matrix
$a:=(a_{i,j})_{i,j\in G}$ is symmetric, and let
\begin{eqnarray*}
w_i(y)&:=&y_i \biggl(2\sum_{j\in S}\bigl(a_{i,j}-2H(x)\bigr)y_j+\sum_{j\in\partial
S}\bigl(a_{i,j}-2H(x)\bigr)y_j \biggr)\\
&\hspace*{3pt}=&o_{|y|\to0}(y_i) =o_{|y|\to0}(y_{\partial S}),
\end{eqnarray*}
using that, for all $j\in\partial S$, $y_j\ge0$.
Finally, we apply in the inequality that $B:=(a_{i,j}-2H(x))_{i,j\in
S}$ is a negative semidefinite matrix by Lemma \ref{lem:1}.

Using that, for all $i\in\partial S$, $N_i(x)<H(x)$ (and $y_i\ge0$),
we deduce that there exists a neighborhood $\Nn(x)$ of $x$ in $\De$
such that, if $x+y\in\Nn(x)$, then $H(x+y)\le H(x)$.\vadjust{\goodbreak}

In order to obtain the required estimate of $H(x+y)-H(x)$, we observe that,
if $z:=(y_i)_{i\in S}$, then, by semidefiniteness of the symmetric matrix~$B$,
%
\begin{equation}
\label{eq:ba}
-\Cst(x,a)|Bz|^2\le(Bz,z)\!=\!\sum_{i,j\in S}\bigl(a_{i,j}-2H(x)\bigr)y_iy_j\le-\Cst(x,a)|Bz|^2.\hspace*{-30pt}
\end{equation}
But
\[
Bz= \biggl(N_i(y)-2H(x)\sum_{i\in S}y_i \biggr)_{i\in S}
=\bigl(N_i(y)+2H(x)y_{\partial S}\bigr)_{i\in S},
\]
where we use that $y_{\partial S}=-y_S$ in the second equality, since
$y\in T\De$. Hence,
%
\begin{equation}
\label{eq:bh}
|Bz|^2=\sum_{i\in S}\bigl(N_i(y)+2H(x)y_{\partial S}\bigr)^2=\sum_{i\in
S}N_i(y)^2+o_{|y|\to0}(y_{\partial S})
\end{equation}
and, if we let
\[
K(y):=\sum_{i\in S}N_i(y)^2+y_{\partial S},
\]
then, by combining identities (\ref{eq:h}), (\ref{eq:ba}) and (\ref
{eq:bh}) [and using that $w_i(y)=o_{|y|\to0}(y_{\partial S})$ for all
$i\in\partial S$], restricting $\Nn(x)$ if necessary,
%
\begin{equation}
\label{eq:hk}
-\Cst(x,a)K(y)\le H(x+y)-H(x)\le-\Cst(x,a)K(y).
\end{equation}
On the other hand, let
\[
L(y):=\sum_{i\in S}\bigl(N_i(x+y)-H(x+y)\bigr)^2+y_{\partial S}.
\]
Then, again by restricting $\Nn(x)$ if necessary,
%
\begin{equation}
\label{eq:jl}
\Cst(x,a)L(y)\le J(x+y)\le\Cst(x,a)L(y),
\end{equation}
where we use again that $N_i(x)<H(x)$ for all $i\in\partial S$.
But
\begin{eqnarray}\label{eq:lk}
L(y)&=&\sum_{i\in S}\bigl[N_i(y)-\bigl(H(x+y)-H(x)\bigr)\bigr]^2+y_{\partial
S}\nonumber\\[-8pt]\\[-8pt]
&=&K(y)+o_{|y|\to0}\bigl(|H(x+y)-H(x)|\bigr).\nonumber
\end{eqnarray}

Combining inequalities (\ref{eq:hk}), (\ref{eq:jl}) and (\ref
{eq:lk}), and further restricting $\Nn(x)$ if necessary, we obtain
inequality (\ref{in:lem4}) as required.

Let us now prove (b). If $q\in\Sc(S(x))$ and $y\in\De$, then
\[
-\sum_{i\in S}q_i[N_i(y)-H(y)]=H(y)-\sum_{i\in S}q_iN_i(y)
\]
and
\[
\sum_{i\in S}q_iN_i(y)=\sum_{i\in G}q_iN_i(y)=\sum_{i\in G}y_iN_i(q)
=H(q)+\sum_{i\in\partial S}y_i[N_i(q)-H(q)],
\]
where we use that $(a_{i,j})_{i,j\in G}$ is symmetric in the second
equality, and that~$q$ is an equilibrium in the third equality.
Therefore,
%
\begin{equation}
\label{var:iq}
\qquad I_q(y)=H(y)-H(q)+\sum_{i\in\partial S}y_i\bigl[2\bigl(N_i(y)-H(y)\bigr)-\bigl(N_i(q)-H(q)\bigr)\bigr].
\end{equation}
If $q, y \in\Nn(x)$, then [by restricting $\Nn(x)$ if necessary]
$x\in\Ec_s$ implies that, for all $i\in\partial S$,
\[
-\Cst(x,a)\le2\bigl(N_i(y)-H(y)\bigr)-\bigl(N_i(q)-H(q)\bigr)\le-\Cst(x,a).
\]
Inequality (\ref{Iq<0}) follows.

\section{Stochastic approximation results for the VRRW}
\label{sec:stoc}

\subsection{The stochastic approximation equation}
\label{stoc-app}

We assume in this section that $G$ is finite. The main idea is to modify the density of occupation measure
\[
x(n)= \biggl(\frac{Z_n(i)}{n+n_0} \biggr)_{i\in G}
\]
into a vector $z(n)$ that takes into account the position of the random
walk, so that the conditional expectation of $z(n+1)-z(n)$ roughly only
depends on $z(n)$ and not on the position $X_n$. This expectation will
actually approximately be $F(z(n))/(n+n_0)$, where $F$ is the map
involved in the ordinary differential equation (\ref{eqdiff}).

For all $x\in\De$, let $M(x)$ be the following matrix of transition
probabilities of the reversible Markov chain:
%
\begin{equation}
\label{markov}
M(x)(i,j)\dvtx \1_{i\sim j} \frac{a_{i,j} x_j}{\sum_{k\sim i} a_{i,k} x_k};
\end{equation}
$M(x(n))$ provides the transition probabilities from the VRRW at time
$n$. Recall that $\pi(x)$ in (\ref{pi}) is the invariant probability
measure for $M(x)$.

Let us denote by $\Gg$ (resp., $\Hh$) the set of functions on $G$
taking values in~$\R$ (resp., in $\R^G$). Let $\1$ be the function
identically equal to $1$.
Let $M(x)$ and~$\Pi(x)$ denote the linear transformations on $\Gg$
defined by
%
\begin{eqnarray}
\label{markov2}
(M(x)f)(i)&:=&\sum_{j\in G}M(x)(i,j)f(j),\\
\label{pix}
\Pi(x)(f)&:= &\biggl(\sum_{i\in G}\pi(x)(i)f(i) \biggr)\1.
\end{eqnarray}
Note that, by a slight abuse of notation, $M(x)$ equally denotes the
Markov chain defined in (\ref{markov}) and its transfer operator in
(\ref{markov2}); $\Pi(x)$ is the linear transformation of $\Gg$ that
maps $f$ to the linear form identically equal to
the mean of $f$ under the invariant probability measure $\pi(x)$.

Any linear transformation $P$ of $\Gg$ [and, in particular, $M(x)$ and
$\Pi(x)$] also defines a linear transformation of $\Hh$: for all
$f=(f_i)_{i\in G}\in\Hh$,
%
\begin{equation}
\label{extend_map}
Pf:=(Pf_i)_{i\in G}.
\end{equation}

Let us now introduce a solution of the Poisson equation for the Markov
chain $M(x)$.
Let us define, for all $t\in\R_+$,
\[
G_t(x):=e^{-t(I-M(x))}=e^{-t}\sum_0^\iy\frac{t^iM(x)^i}{i!},
\]
which is the Markov operator of the continuous time Markov chain
associated with $M(x)$. For all $x\in \mathit{Int}(\De)$, $M(x)$ is
indecomposable so that $G_t(x)$ converges toward $\Pi(x)$ at an
exponential rate, hence,
\[
Q(x):=\int_0^\iy\bigl(G_t(x)-\Pi(x)\bigr) \,dt
\]
is well defined. Note that
\[
Q(x)\1=0,
\]
and that $Q(x)$ is the solution of the Poisson equation
%
\begin{equation}
\label{poisson}
\bigl(I-M(x)\bigr)Q(x)=Q(x)\bigl(I-M(x)\bigr)=I-\Pi(x),
\end{equation}
using that $M(x)\Pi(x)f=\Pi(x)f=\Pi(x)M(x)f$ for all $f\in\Gg$ (or
$f\in\Hh$).

Let us now expand $x(n+1)-x(n)$, using (\ref{poisson}). Let
$(e_i)_{i\in G}$ be the canonical basis of $\R^G$, that is, $e_i:=(\1
_{j=i})_{j\in G}$ for all $i\in G$. Let $\iota\in\Hh$ be defined by
\begin{eqnarray*}
\iota\dvtx   G &\longrightarrow&\R^G,\\
  i  &\longmapsto& e_i.
\end{eqnarray*}

First note that, for all $x\in\De$,
$\Pi(x)\iota=\pi(x)\1$ since, for all $j\in G$,
\[
\Pi(x)\iota(j)=((\Pi(x)\iota_k)(j))_{k\in G}
=((\pi(x)(k)\1)(j))_{k\in G}=\pi(x).
\]

Therefore,
\begin{eqnarray*}
(n+n_0+1)\bigl(x(n+1)-x(n)\bigr) &=&\bigl(\1_{X_{n+1}=i}-x(n)_i\bigr)_{i\in G}
=\iota(X_{n+1})-x(n)\\
&=&\iota(X_{n+1})-\pi(x(n))+F(x(n))\\
&=&[I-\Pi(x(n))]\iota(X_{n+1})+F(x(n)),
\end{eqnarray*}
where $F$ is the function defined in  (\ref{replicator}).

Now,
\begin{eqnarray}
\label{poissondecomp}
\frac{[I-\Pi(x(n))]\iota(X_{n+1})}{n+n_0+1}
&=&\frac{(Q(x(n))-M(x(n)Q(x(n)))\iota(X_{n+1})}{n+n_0+1}\nonumber\\[-8pt]\\[-8pt]
&=&\epsilon_{n+1}+\eta_{n+1}+r_{n+1,1}+r_{n+1,2},\nonumber
\end{eqnarray}
where
\begin{eqnarray*}
\epsilon_{n+1}&:=&
\frac{Q(x(n))\iota(X_{n+1})-M(x(n))Q(x(n))\iota(X_n)}{n+n_0+1},
\\[-1pt]
r_{n+1,1}&:=&\biggl(\frac{1}{n+n_0+1}-\frac{1}{n+n_0}
\biggr)M(x(n))Q(x(n))\iota(X_n)\\[-1pt]
&\hspace*{3pt}=&-\frac{M(x(n))Q(x(n))\iota(X_n)}{(n+n_0)(n+n_0+1)},\\[-1pt]
\eta_{n+1}&:=&\frac{M(x(n))Q(x(n))\iota(X_n)}{n+n_0}-
\frac{M(x(n+1))Q(x(n+1))\iota(X_{n+1})}{n+n_0+1},\\[-1pt]
r_{n+1,2}&:=&
\frac{[M(x(n+1))Q(x(n+1))-M(x(n))Q(x(n))]\iota(X_{n+1})}{n+n_0+1}.
\end{eqnarray*}

Let, for all $n\in\N$,
%
\begin{equation}
\label{defzn}
z(n):=x(n)+\frac{M(x(n))Q(x(n))\iota(X_n)}{n+n_0}
\end{equation}
and
\begin{eqnarray*}
r_{n+1,3}&:=&\frac{1}{n+n_0+1}\frac{F(x(n))-F(z(n))}{H(x(n))},\\[-1pt]
r_{n+1}&:=&r_{n+1,1}+r_{n+1,2}+r_{n+1,3}.
\end{eqnarray*}

Then, for all $n\in\N$, it follows from equation (\ref
{poissondecomp}) that
%
\begin{equation}
\label{eq stoc-approx}
z(n+1)=z(n)+\frac{1}{n+n_0+1}\frac{F(z(n))}{H(x(n))}+\epsilon_{n+1}+r_{n+1}.
\end{equation}
Note that $\Es(\epsilon_{n+1} | \F_n)=0$, since
\[
\Es(Q(x(n))\iota
(X_{n+1}) | \F_n)=M(x(n))Q(x(n))\iota(X_n);
\]
also observe that
\[
\sum_{i\in G}z(n)_i=\sum_{i\in G}x(n)_i+\frac{(M(x(n))Q(x(n))\1
)(X_n)}{n+n_0}=1.
\]
We provide in Section \ref{as-trap} estimates of the conditional
variance of $\epsilon_{n+1}$ and of $r_{n+1}$, which will be
sufficient to prove localization of the vertex-reinforced random walk
with positive probability.\vspace*{-2pt}

\subsection{Estimates on the underlying Markov chain $M(x)$}
\label{as-trap}
For convenience we assume here that $G=S\cup\partial S$, where
$(S,\sim)$ is finite, connected and not a~singleton unless it is a
loop. Let $\oa:=\max_{i,j\in G, i\sim j}a_{i,j}$,\vspace*{2pt}
$\ua:=\min_{i,j\in G, i\sim j}a_{i,j}$.

Let us first introduce some general notation on Markov chains. Let $K$
be a reversible Markov chain on the graph $(G,\sim)$, with invariant
measure $\mu$. Let $\langle\cdot,\cdot\rangle_\mu$ be the scalar product defined by, for
all $f, g \in\Gg$,
\[
\langle f,g\rangle_\mu:=\sum_{x\in G}f(x)g(x)\mu(x).\vadjust{\goodbreak}
\]
On $\Gg$, we define the $\ell^p(\mu)$ norm, $1\le p<\iy$ by
\[
\|f\|_{\ell^p(\mu)}:= \biggl(\sum_{x\in G}|f(x)|^p\mu(x) \biggr)^{1/p},\vspace*{-1.5pt}
\]
and the infinity norm
\[
\|f\|_\iy:=\max_{x\in G} |f(x)|.\vspace*{-1.5pt}
\]
We also define the infinity norm on $\Hh$: if $f=(f_i)_{i\in G}\in\Hh$,
%
\begin{equation}
\label{extend_norm}
\|f\|_\iy=\max_{i\in G}\|f_i\|_\iy=\max_{i,x\in G}|f_i(x)|.\vspace*{-1.5pt}
\end{equation}

Let $\Es_\mu$ be the expectation operator
\[
\Es_\mu f:=\sum_{x\in G} f(x)\mu(x)=\langle f,\1\rangle_\mu,\vspace*{-1.5pt}
\]
where $\1$ is the constant function equal to $1$.

We let
$\E_K$ be the Dirichlet form of $K$,
\[
\E_K(f,g)=\langle(I-K)f,g\rangle_\mu,\vspace*{-1.5pt}
\]
and let $\Var_\mu$ be the variance operator,
\[
\Var_\mu(f):=\|f-\Es_\mu f\|_{\ell^2(\mu)}^2=\|f\|_{\ell^2(\mu
)}^2-(\Es_\mu f)^2.\vspace*{-1.5pt}
\]

Simple calculations yield that
\[
\E_K(f,f)=\frac{1}{2}\sum_{i\sim j}\bigl(f(i)-f(j)\bigr)^2K(i,j)\mu(i),\vspace*{-1.5pt}
\]
and
\[
\Var_\mu(f)=\frac{1}{2}\sum_{i, j\in G}\bigl(f(i)-f(j)\bigr)^2\mu(i)\mu(j).\vspace*{-1.5pt}
\]

Let $\lambda(K)$ be the spectral gap of the Markov chain $K$,
\[
\lambda(K):=\min\biggl\{\frac{\E_K(f,f)}{\Var_\mu(f)}\mbox{ s.t. }
\Var_\mu(f)\not=0 \biggr\}.\vspace*{-1.5pt}
\]

The following Lemma \ref{sp_gap} states that the spectral gap of the
Markov chain~$M(x)$ is lower bounded on $\La_\al$ [defined in (\ref
{eq:lalph})].
\begin{lemma}
\label{sp_gap}
For all $x\in\La_\al$, $\lambda(M(x))\ge\Cst(\al,a,|G|)$.
\end{lemma}
\begin{pf}
Let $M:=M(x)$ and $\pi:=\pi(x)$ for simplicity. Let us first observe
that, for all $i\in G$, $j\in S$ such that $i\sim j$,
\begin{eqnarray}
\label{pi_M}
M(i,j)&\ge&\ua x_j/\oa\ge\al\ua/\oa\quad \mbox{and }\nonumber\\[-9pt]\\[-9pt]
M(i,j)\pi(i)&=&\pi
(j)M(j,i)\ge\ua\al^2\1_{i\in S}/\oa,\nonumber\vadjust{\goodbreak}
\end{eqnarray}
where the second inequality comes from
\[
M(i,j)\pi(i)=\frac{a_{i,j}x_j}{N_i(x)}\frac{x_iN_i(x)}{H(x)}
=\frac{a_{i,j}x_ix_j}{H(x)}\ge\frac{\ua\al^2}{\oa}\1_{i\in S}.
\]

Now, by connectedness of $(S,\sim)$, for all $i, j \in G$, there
exists $l\le|G|$ and a path
$(n_k)_{1\le k\le l}\in G\times S^{l-2}\times G$ such that $i=n_1$,
$j=n_l$, $n_k\sim n_{k+1}$ for all $k\in\{1,\ldots,l-1\}$.

Hence, for all $k\in\{1,\ldots,l\}$, using inequalities (\ref{pi_M}),
\begin{eqnarray*}
&&\pi(i)\pi(j)\bigl(f(i)-f(j)\bigr)^2\\
&&\qquad \le l\pi(i)\pi(j)\sum_{k\in\{1,\ldots,l-1\}}\bigl(f(n_k)-f(n_{k+1})\bigr)^2\\
&&\qquad \le l\pi(i)\bigl(f(i)-f(n_2)\bigr)^2+
l\pi(j)\bigl(f(j)-f(n_{l-1})\bigr)^2\\
&&\qquad \quad {}+
l\sum_{k\in\{2,\ldots,l-2\}}\bigl(f(n_k)-f(n_{k+1})\bigr)^2\\
&&\qquad \le\frac{\oa l}{\ua\al}\bigl[M(i,n_2)\pi(i)\bigl(f(i)-f(n_2)\bigr)^2\\
&&\qquad \quad \hphantom{\frac{\oa l}{\ua\al}\bigl[}
{}+M(j,n_{l-1})\pi(j)\bigl(f(j)-f(n_{l-1})\bigr)^2\bigr]\\
&&\qquad \quad {} +\frac{\oa l}{\ua\al^2}\sum_{k\in\{2,\ldots,l-2\}
}\bigl(f(n_k)-f(n_{k+1})\bigr)^2M(n_k,n_{k+1})\pi(n_k)\\
&&\qquad \le\frac{\oa l}{\ua\al^2}\sum_{k\in\{1,\ldots,l-1\}
}\bigl(f(n_k)-f(n_{k+1})\bigr)^2M(n_k,n_{k+1})\pi(n_k)\\
&&\qquad \le\frac{2\oa|G|}{\ua\al^2}\E_M(f,f).
\end{eqnarray*}
Therefore,
\[
\Var_\pi(f)=\frac{1}{2}\sum_{i, j\in G}\pi(i)\pi
(j)\bigl(f(i)-f(j)\bigr)^2\le\frac{\oa|G|^3}{\ua\al^2}\E_M(f,f).
\]
\upqed
\end{pf}

Lemma \ref{bound Q(x)} provides upper bounds on the
norms of $Q(x)$, $M(x)Q(x)$ and their partial derivatives on $\La_\al
$, which will be needed in the estimates of $r_{n+1}$ and of the
conditional variance of $\epsilon_{n+1}$ in Lemma \ref{bound_rn}.

The norm on linear transformations of $\Gg$ will be the infinity norm
\[
\|A\|_\iy:=\sup_{f\in\Gg, f\not=0}\frac{\|Af\|_\iy}{\|f\|_\iy}.
\]
Note that, for any linear transformation $A$ of $\Gg$, the
corresponding linear transformation of $\Hh$ (still called $A$)\vadjust{\goodbreak}
defined in (\ref{extend_map}) still has the same infinity norm [the $\|
\cdot\|_\iy$ on $\Hh$ is defined by (\ref{extend_norm})],
\[
\|A\|_\iy=\sup_{f\in\Hh, f\not=0}\frac{\|Af\|_\iy}{\|f\|_\iy}.\vspace*{-3pt}
\]

\begin{lemma}
\label{bound Q(x)}
For all $x\in\La_\al$, $i, j \in G$, $f\in\Gg$:
\begin{eqnarray*}
&&\mathrm{(a)}\quad M(x)(i,j)\le\biggl(\frac{\oa}{\ua}
\biggr)^2\frac{\pi(x)(j)}{\al^2},\\[-1.5pt]
&&\mathrm{(b)} \quad \|Q(x)f\|_{\ell^2(\pi(x))}\le\frac{\sqrt{\Var_{\pi
(x)}(f)}}{\lambda(M(x))}
\le\frac{\|f\|_{\ell^2(\pi(x))}}{\lambda(M(x))},\\[-1.5pt]
&&\mathrm{(c)} \quad \|Q(x)\|_\iy\le\Cst(\al,a,|G|),\qquad
\|M(x)Q(x)\|_\iy\le\Cst(\al,a,|G|),\\[-1.5pt]
&&\mathrm{(d)} \quad \biggl\|\frac{\partial Q(x)}{\partial x_i}\biggr\|_\iy\le\Cst(\al,a,|G|),\qquad
\biggl\|\frac{\partial(M(x)Q(x))}{\partial x_i}\biggr\|_\iy\le\Cst(\al,a,|G|).
\end{eqnarray*}
\end{lemma}
\begin{pf}
Let $M:=M(x)$, $Q:=Q(x)$, $\pi:=\pi(x)$, $\lambda:=\lambda(M(x))$
for simplicity.

Inequality (a) is obvious: for all $j\in G$,
\[
M(i,j)=\frac{a_{i,j}x_j}{N_i(x)}=\frac{x_jN_j(x)}{H(x)}\frac
{a_{i,j}H(x)}{N_i(x)N_j(x)}
\le\biggl(\frac{\oa}{\ua} \biggr)^2\frac{\pi(j)}{\al^2}.
\]

Let us now prove (b). For all $f\in\Gg$,
\[
\|G_tf-\pi(f)\|_{\ell^2(\pi)}^2\le e^{-2\lambda t}\Var_\pi(f),
\]
by definition of the spectral gap (see, e.g., Lemma 2.1.4, \cite
{saloff}), so that
%
\begin{eqnarray}\label{l2norm}
\nonumber
\|Q(x)f\|_{\ell^2(\pi)}&\le&\biggl\|\int_0^\iy\bigl(G_t(x)f-\Pi(x)f\bigr) \,dt \biggr\|
_{\ell^2(\pi)}\\[-1.5pt]
&\le&\int_0^\iy\bigl\|\bigl(G_t(x)f-\Pi(x)f\bigr) \bigr\|_{\ell^2(\pi)} \,dt\\[-1.5pt]
&\le&\sqrt{\Var_\pi(f)}\int_0^\iy e^{-\lambda t} \,dt=\frac{\sqrt
{\Var_\pi(f)}}{\lambda}
\le\frac{\|f\|_{\ell^2(\pi)}}{\lambda}.\nonumber
\end{eqnarray}
Inequality (c) translates this upper bound of the $\ell^2(\pi
)\to\ell^2(\pi)$-norm of $Q(x)$ into one involving the infinity norm
for $MQ$, using (a):
\begin{eqnarray*}
|MQf(i)|&=& \biggl|\sum_{j\in G} M(i,j)Qf(j) \biggr|\\[-1.5pt]
&\le&\frac{1}{\al^2} \biggl(\frac{\oa}{\ua} \biggr)^2\sum_{j\in G}\pi(j)|Qf(j)|=
\biggl(\frac{\oa}{\ua} \biggr)^2\frac{\|Qf\|_{\ell^1(\pi)}}{\al^2}\\[-1.5pt]
&\le& \biggl(\frac{\oa}{\ua} \biggr)^2\frac{\|Qf\|_{\ell^2(\pi)}}{\al^2}\le
\biggl(\frac{\oa}{\ua}
\biggr)^2\frac{\|f\|_{\ell^2(\pi)}}{\lambda\al^2}.\vadjust{\goodbreak}
\end{eqnarray*}

Hence, using Lemma \ref{sp_gap},
\[
\|MQf\|_\iy\le\biggl(\frac{\oa}{\ua} \biggr)^2\frac{\|f\|_{\ell^2(\pi
)}}{\lambda\al^2}\le\biggl(\frac{\oa}{\ua} \biggr)^2\frac{\|f\|_{\iy
}}{\lambda\al^2}
\le\Cst(\al,a,|G|)\|f\|_{\iy}.
\]
Then the same upper bound for $\|Q(x)f\|_\iy$ follows from the Poisson
equation~(\ref{poisson}):
\[
Q(x)=M(x)Q(x)+I-\Pi(x).
\]

Let us now prove (d). Given $i\in G$, let us take the derivative
of the Poisson equation $Q(x)(I-M(x))=I-\Pi(x)$ with respect to $x_i$:
\[
\frac{\partial Q(x)}{\partial x_i} \bigl(I-M(x)\bigr)=Q(x)\frac{\partial
M(x)}{\partial x_i}
-\frac{\partial\Pi(x)}{\partial x_i}.
\]
This equality, multiplied on the right by $Q(x)$, yields, using now the
Poisson equation $(I-M(x))Q(x)=I-\Pi(x)$,
%
\begin{equation}
\label{deriv_poisson}
\frac{\partial Q(x)}{\partial x_i}=\frac{\partial Q(x)}{\partial
x_i}\bigl(I-\Pi(x)\bigr)
= \biggl(Q(x)\frac{\partial M(x)}{\partial x_i}-\frac{\partial\Pi
(x)}{\partial x_i} \biggr)Q(x),
\end{equation}
where we use that, for all $f\in\Gg$,
\[
\frac{\partial Q(x)}{\partial x_i}\Pi(x)f=\langle f,\1\rangle_{\pi(x)}\frac
{\partial Q(x)}{\partial x_i}\1=0,
\]
since $Q(x)\1=0$ for all $x\in\De$.

Equality (\ref{deriv_poisson}) implies the required upper bound of
$\|\frac{\partial Q(x)}{\partial x_i}\|_\iy$. Indeed, the following
estimates hold: for all $i, j, k \in G$, $j\sim k$,
\begin{eqnarray*}
\biggl|\frac{\partial[M(x)(j,k)]}{\partial x_i} \biggr|
&=& \biggl|\frac{\partial}{\partial x_i} \biggl(\frac{a_{j,k}x_k}{N_j(x)} \biggr)
\biggr|\\
&=& \biggl|\frac{\partial x_k}{\partial x_i}\frac{a_{j,k}}{N_j(x)}-\frac
{a_{j,k}x_k}{N_j(x)^2}\frac{\partial N_j(x)}{\partial x_i} \biggr|\\
&\le&\frac{2\oa}{N_j(x)}\le\frac{2\oa}{\ua\al},
\end{eqnarray*}
where we use that $a_{j,k}x_k\le N_j(x)$ and $\partial N_j/\partial
x_i(x)=a_{j,i}$, and that there exists $l\in S$ with $l\sim j$, given
the assumptions on $S$. Also,
\begin{eqnarray*}
\biggl|\frac{\partial\pi(x)(j)}{\partial x_i} \biggr|&=&
\biggl|\frac{\partial}{\partial x_i} \biggl(\frac{x_jN_j(x)}{H(x)} \biggr)
\biggr|\\
&=& \biggl|\frac{\partial(x_jN_j(x))}{\partial x_i}\frac{1}{H(x)}-\frac
{x_jN_j(x)}{H(x)^2}\frac{\partial H(x)}{\partial x_i} \biggr|\\
&\le&\frac{4\oa}{H(x)}\le\frac{4\oa}{\ua\al^2},
\end{eqnarray*}
where we note that $|\frac{\partial H(x)}{\partial
x_i}|\,{=}\,2N_i(x)\,{\le}\,2\oa$. The upper bound of $\|\frac{\partial(M(x)Q(x))}{\partial
x_i}\|_\iy$ follows directly.
\end{pf}

\subsection{\texorpdfstring{Proof of Lemma \protect\ref{bound_rn}}{Proof of Lemma 4}}
\label{sec:bound rn}
The estimates (a) and (d) readily follow from the
definitions of $\epsilon_{n+1}$ and $z(n)$, and from Lemma \ref{bound
Q(x)}(c).

Let $M:=M(x(n))$, $Q:=Q(x(n))$, $\pi:=\pi(x(n))$, $\lambda:=\lambda
(M(x(n)))$ for simplicity.
Let us prove (b):
\begin{eqnarray*}
(n+n_0)^2\Es((\epsilon_{n+1})_i^2 | \F_n)
&\le&\Es([Qe_i(X_{n+1})]^2 | \F_n)\\
&=&\sum_{j\sim X_n} M(X_n,j)[Qe_i(j)]^2\\
&\le&\frac{1}{\al^2} \biggl(\frac{\oa}{\ua} \biggr)^2\sum_{j\in G}
\pi(j)[Qe_i(j)]^2\\
&=& \biggl(\frac{\oa}{\ua} \biggr)^2\frac{1}{\al^2}\|Qe_i\|_{\ell^2(\pi
(x(n))}^2\\
&\le&\Cst(\al,a,|G|)\|e_i\|_{\ell^2(\pi(x(n))}^2
\le\Cst(\al,a,|G|)x(n)_i,
\end{eqnarray*}
where we use Lemma \ref{bound Q(x)}(a) and (b),
respectively, in the second and in the third inequality.

In order to prove (c), let us first upper bound $\|r_{n+1,1}\|
_\iy$ using Lemma~\ref{bound Q(x)}(c):
\[
\|r_{n+1,1}\|_\iy
\le\frac{\|M(x(n))Q(x(n))\iota(X_n)\|_\iy}{(n+n_0)^2}
\le\frac{\Cst(\al,a,|G|)}{(n+n_0)^2}.
\]

Let us now bound $\|r_{n+1,2}\|_\iy$:
\begin{eqnarray*}
(n+n_0)\|r_{n+1,2}\|_\iy
&\le&\sup_{\theta\in[0,1]} \biggl\|\frac{\partial(MQ)(\theta
x(n)+(1-\theta)x(n+1))}{\partial\theta} \biggr\|_\iy\\
&\le&\sum_{i\in G} \bigl|\bigl(x(n+1)-x(n)\bigr)_i\bigr|\\
&&\hphantom{\sum_{i\in G}}
{}\times\sup_{i\in G, \theta\in[0,1]}
\biggl\|\frac{\partial(MQ)(\theta x(n)+(1-\theta)x(n+1))}{\partial x_i} \biggr\|
_\iy\\
&\le&\frac{\Cst(\al,a,|G|)}{n+n_0},
\end{eqnarray*}
where we use Lemma \ref{bound Q(x)}(d) in the last inequality.

It remains to upper bound $\|r_{n+1,3}\|_\iy$. First observe that, for
all~\mbox{$y\,{=}\,(y_i)_{i\in G}$}, $z=(z_i)_{i\in G} \in\De$, $i\in G$,
\[
|F_i(z)-F_i(y)|\le\sum_{j\in G}|z_j-y_j|\sup_{k\in G, x\in\De}
\biggl|\frac{\partial F_i(x)}{\partial x_k} \biggr|\le2\oa\sum_{i\in G}|z_i-y_i|,
\]
where we use the explicit computations of $\partial F_i/\partial x_j$
in the proof of Lemma~\ref{lem:1}. Hence,
\[
\|F(z)-F(y)\|_\iy\le2\oa|G|\|z-y\|_\iy,\vadjust{\goodbreak}
\]
which implies
\[
\|r_{n+1,3}\|_\iy\le
\frac{1}{n+n_0}\frac{|G|}{\ua}2\oa|G|\|x(n)-z(n)\|_\iy
\le\frac{\Cst(\al,a,|G|)}{(n+n_0)^2},
\]
where we use that, by inequality (\ref{lbound_H}), $H(x)\ge\ua/|G|$
for all $x\in\De$.

\subsection{\texorpdfstring{Proof of Lemma \protect\ref{lem:evolent} and inclusions (\protect\ref{eqVqiy2})}
{Proof of Lemma 5 and inclusions (15)}}
\label{sec:evolent}
Let us first prove inclusions~(\ref{eqVqiy2}). If we let $g\dvtx \R
_+\setminus\{0\}\longrightarrow\R_+$ be the function defined by
$g(u):=u-\log(u+1)$, nonnegative by concavity of the log function,
then, for all $y\in\De$ such that $y_i>0$ for all $i\in S$,
%
\begin{equation}\label{eqVqiy1}
\qquad V_q(y)=-\sum_{i\in S}q_i\log\biggl(1+\frac{y_i-q_i}{q_i} \biggr)+2y_{\partial S}
=\sum_{i\in S}q_i g \biggl(\frac{y_i-q_i}{q_i} \biggr)+3y_{\partial S},
\end{equation}
which implies the inclusions.

Let us now prove Lemma \ref{lem:evolent}; let, for all $n\in\N$,
\[
\ze_{n+1}:= \biggl(\frac{(\epsilon_{n+1})_i}{z(n)_i}\1_{i\in S} \biggr)_{i\in G},
\]
with the convention that $\ze_{n+1}=0$ if $z(n)_i=0$ for some $i\in S$.
Fix $\epsilon>0$ such that $B_{V_q}(2\epsilon)\subset\La_\al$ for
some $\al=\Cst(q)>0$, and assume $x(n)\in B_{V_q}(\epsilon)$ for
some $n\ge n_1$. Thus, $\|z(n)-x(n)\|_\iy\le\Cst(q,a)/(n+n_0)$ by
Lemma \ref{bound_rn}(d); we assume in the rest of the proof
that $\epsilon<\Cst(q)$ and $n_0\ge\Cst(q,a)$ so that, using (\ref
{eqVqiy1}), $z(n)\in B_{V_q}(2\epsilon)\subset\La_\al$.

Note that $\|x(n)-x(n+1)\|_\iy\le(n+n_0)^{-1}$, which implies, using
Lemma \ref{bound_rn}, that $\|z(n)-z(n+1)\|_\iy\le\Cst
(q,a)(n+n_0)^{-1}$. Hence, using that $z(n)\in\La_\al$,
\begin{eqnarray*}
V_q\bigl(z(n+1)\bigr)-V_q(z(n))
&=&-\sum_{i\in S}q_i\log\biggl(\frac{z(n+1)_i}{z(n)_i}
\biggr)+2[z(n+1)_{\partial S}-z(n)_{\partial S}]\\
&=&-\sum_{i\in S}q_i\frac{z(n+1)_i-z(n)_i}{z(n)_i}+2[z(n+1)_{\partial
S}-z(n)_{\partial S}]\\
&&{}+\Box\biggl(\frac{\Cst(q,a)}{(n+n_0)^2} \biggr),
\end{eqnarray*}
where we again make use of notation $u=\Box(v)\iff|u|\le v$ from
Section~\ref{sec:notation}.

Hence, using identity (\ref{eq stoc-approx}) and Lemma \ref
{bound_rn}(c)--(d), we obtain subsequently [recall that $I_q(\cdot)$ is defined
in (\ref{eq:iq})]
\begin{eqnarray*}
&&V_q\bigl(z(n+1)\bigr)-V_q(z(n))\\
&&\qquad =\frac{1}{n+n_0+1}\frac
{I_q(z(n))}{H(x(n))}-(q,\ze_{n+1})+2(\epsilon_{n+1})_{\partial
S}+\Box\biggl(\frac{\Cst(q,a)}{(n+n_0)^2} \biggr)\\
&&\qquad =\frac{1}{n+n_0+1}\frac{I_q(x(n))}{H(x(n))}-(q,\ze
_{n+1})+2(\epsilon_{n+1})_{\partial S}+\Box\biggl(\frac{\Cst
(q,a)}{(n+n_0)^2} \biggr).
\end{eqnarray*}

\subsection{\texorpdfstring{Proof of Lemma \protect\ref{lem:stoc h}}{Proof of Lemma 9}}
\label{sec:stoc h}
Using identities (\ref{eq:h0}) and (\ref{eq stoc-approx}) [recall
that $J$ is defined in (\ref{defj})],
\begin{eqnarray*}
H\bigl(z(n+1)\bigr)-H(z(n))&=&2\sum_{i\in
G}N_i(z(n))\cdot\bigl(z(n+1)-z(n)\bigr)_i\\
&&{}+H\bigl(z(n+1)-z(n)\bigr)\\
&=&
\frac{1}{n+n_0+1}\frac{J(z(n))}{H(x(n))}+\xi_{n+1}+s_{n+1},
\end{eqnarray*}
where
\begin{eqnarray*}
\xi_{n+1}&:=&2\sum_{i\in G}N_i(z(n))(\epsilon_{n+1})_i,\\
s_{n+1}&:=&2\sum_{i\in G}N_i(z(n))(r_{n+1})_i+H\bigl(z(n+1)-z(n)\bigr).
\end{eqnarray*}

Let $\al>0$, and assume $x(n)\in\La_\al$. Inequalities
(\ref{nih}) and (\ref{pi}) of our lemma follow from Lemma~\ref{bound_rn}(a)--(c), and from $\|
z(n+1)-z(n)\|_\iy\le\Cst(\al,a,|G|)/\allowbreak(n+n_0)$ (see, e.g.,
the beginning of the proof of Lemma \ref{lem:evolent}).

\section{Asymptotic results for the VRRW}
\label{sec:asymp}

\subsection{\texorpdfstring{Proof of Lemma \protect\ref{lem:posprob}}{Proof of Lemma 7}}
\label{sec:enttrap}
Fix $\epsilon>0$ such that $B_{V_x}(\epsilon)\subset\La_\al$ for
some $\al>0$ depending on $x$, and assume $x(n)\in B_{V_x}(\epsilon
/2)$ for some $n\ge n_1$.

Let $(\ze_k)_{k\ge2}$ be defined as in Section \ref{sec:evolent},
and let us define the martingales $(A_k)_{k\ge n }$, $(B_k)_{k\ge n }$
and $(\ka_k)_{k\ge n }$ by
\begin{eqnarray*}
A_k&:=&\sum_{j=n+1}^k \ze_j\1_{\{V_x(x(j-1))\le\epsilon\}}, \qquad
B_k:=\sum_{j=n+1}^k (\epsilon_j)_{\partial S}\1_{\{V_x(x(j-1))\le
\epsilon\}},\\
\ka_k&:=&-(q,A_k)+2B_k,
\end{eqnarray*}
with the convention that $A_n:=0$ and $B_n=\ka_n:=0$. Using Lemma \ref
{bound_rn}(a), it follows from Doob's convergence theorem that
$(A_k)_{k\ge n }$, $(B_k)_{k\ge n }$ and $(\ka_k)_{k\ge n }$ converge
a.s. and in $\Lc^2$.

Let us briefly outline the proof: we first show that, on an event of
large probability $\Ep$, where $\ka_k$, $k\ge n$, remains small,
$x(k)$ remains in the neighborhood of $x$ and the stochastic
approximation \eqref{evolent} remains valid. This implies, together
with \eqref{Iq<0}, the existence of a subsequence $j_k$ such that
$(x(j_k))_{k\ge0}$ converges to a random $r\in\E_s$ [see \eqref
{eq:subseq}]. Using the linearity of the martingale part of \eqref
{evolent} in $\ze$ and $\epsilon$, we can conclude from the a.s.
convergence of $(A_k)_{k\ge n}$ and $(B_k)_{k\ge n}$ that $x(k)\des
_{k\to\iy}r$ a.s. [see \eqref{eq:vrz} and \eqref{eq:conc}].

The upper bound $|\ka_k-\ka_{k-1}|\le\G/(k+n_0)$ a.s., for some $\G
:=\Cst(x,a)$, implies that, for all $k\ge n+1$ and $\theta\in\R$,
\[
\Es\bigl(\exp\bigl(\theta(\ka_k-\ka_{k-1})\bigr) | \F_{k-1}\bigr)\le\exp\biggl(\frac{\G
^2}{2}\frac{\theta^2}{(k+n_0)^2} \biggr).\vadjust{\goodbreak}
\]
On the other hand, $(\exp(\theta\ka_k))_{k\ge n}$ is a submartingale
since $(\ka_k)_{k\ge n }$ is a~martingale, so that Doob's
submartingale inequality implies, for all $\theta>0$,
\begin{eqnarray*}
\Pb\Bigl(\sup_{k\ge n}\ka_k\ge c | \F_n \Bigr)&=&
\Pb\Bigl(\sup_{k\ge n}e^{\theta\ka_k}\ge e^{\theta c} | \F_n \Bigr)
\le e^{-\theta c}\Es(e^{\theta\ka_\iy} | \F_n)\\
&\le& \exp\biggl(-\theta c+\frac{\theta^2\G^2}{2(n+n_0)} \biggr).
\end{eqnarray*}
Choosing $\theta:=c(n+n_0)/\G^2$ yields
%
\begin{equation}
\label{in_out}
\Pb\Bigl(\sup_{k\ge n}\ka_k\ge c \big| \F_n \Bigr)\le\exp\biggl(-\frac{c^2}{2\G
^2}(n+n_0) \biggr).
\end{equation}

Let
\[
\Ep:= \biggl\{\sup_{k\ge n}\ka_k<\frac{\epsilon}{12} \biggr\};
\]
inequality (\ref{in_out}) implies that
\[
\Pb(\Ep| \F_n)\ge1-\exp\bigl(-\epsilon^2\Cst(x,a)(n+n_0)\bigr).
\]

Now assume that $\Ep$ holds, and let $T$ be the stopping time
\[
T:=\inf\{k\ge n\mbox{ s.t. }V_x(z(k))\ge2\epsilon/3\}.
\]
Note that, using Lemma \ref{bound_rn}(d), if $n\ge\Cst(x,a)$,
then for all $k\in[n,T)$,\break \mbox{$V_x(x(k))<\epsilon$}.
We upper bound $V_x(x(T))-V_x(x(k))$ by adding up identi\-ty~(\ref
{evolent}) in Lemma~\ref{lem:evolent} with $q:=x$, from time $n$ to
$T-1$: this yields, together with Lemma \ref{lem:hestimates}, that
$V_x(z(T))<2\epsilon/3$ if
$T<\iy$, if we assume $n\ge n_1:=\Cst(x,a)$ large enough and
$\epsilon<\epsilon_0:=\Cst(x,a)$ small enough.

Therefore, $V_x(x(k))<\epsilon$ for all $k\ge n$. Using again identity
(\ref{evolent}) [and Lemma~\ref{lem:hestimates}(b)], we obtain
subsequently that
\[
\liminf_{k\to\iy} [H(x)-H(x(k))+x(k)_{\partial S} ]=0\qquad \mbox{a.s.}
\]
since, otherwise, the convergence of $(\ka_k)$ as $k\to\iy$ would
imply\break
$\lim_{k\to\iy}V_x(z(k))=\lim_{k\to\iy}V_x(x(k))=-\iy$,
which is in contradiction with $V_x(x(k))\ge0$.

Hence, there exists a (random) increasing sequence $(j_k)_{k\ge0}$
such that
%
\begin{equation}
\label{eq:subseq}
\lim_{k\to\iy}H(x(j_k))=H(x),\qquad  \lim_{k\to\iy}x(j_k)_{\partial S}=0.
\end{equation}

Let $r$ be an accumulation point of $(x(j_k))_{k\ge0}$. Then
$H(r)=H(x)$ and $r_{\partial S}=0$.

Note that $V_x(r)=\lim_{k\to\iy}V_x(z(j_k))\le\epsilon$. By
possibly choosing a smaller $\epsilon_0:=\Cst(x,a)$, we obtain by
Lemma \ref{lem:hestimates} that $r$ is an equilibrium, and by Lemma
\ref{lem:1} that it is strictly stable.

Let, for all $j\in\N$,
\[
\La_j:= \biggl\{\sup_{k\ge j}|A_k-A_j|<\frac{\epsilon}{24} \biggr\}\cap\biggl\{
\sup_{k\ge j}|B_k-B_j|<\frac{\epsilon}{24} \biggr\}.\vadjust{\goodbreak}
\]

There exists a.s. $j\in\N$ such that $\La_j$ holds; let $l_0$ be such a $j$ ($l_0$ is random, and is not a stopping
time).

Let $k\in\N$ be such that $j_k\ge l_0$ and $V_r(z(j_k))<\epsilon/2$.
Then Lemma \ref{lem:evolent} applies to $r\in\Sc\cap\Ec_s$ and a
similar argument as previously shows that, for all $j'\ge j\ge j_k$,
$V_r(x(j))\le\epsilon$ and
%
\begin{equation}
\label{eq:vrz}
V_r(z(j'))\le V_r(z(j))+\sup_{k\ge j}|A_k-A_j|+2\sup_{k\ge
j}|B_k-B_j|+\frac{\Cst(q,a)}{j+n_0},\hspace*{-20pt}
\end{equation}
if $n_1:=\Cst(x,a)$ was chosen sufficiently large.

Now, $\liminf_{j\to\iy}V_r(z(j))=0$ and
%
\begin{equation}
\label{eq:conc}
\lim_{j\to\iy}\sup_{k\ge j}|A_k-A_j|=\lim_{j\to\iy}\sup_{k\ge
j}|B_k-B_j|=\lim_{j\to\iy}\frac{\Cst(q)}{j+n_0}=0,
\end{equation}
hence, $\lim_{j\to\iy}V_r(x(j))=0$ which implies $\lim_{j\to\iy
}x(j)=r$ and completes the proof.

\subsection{\texorpdfstring{Proof of Lemma \protect\ref{lem:speed}}{Proof of Lemma 8}}
\label{sec:speed}
Let us start with an estimate of the rate of convergence of $H(z(n))$
to $H(x)$.
Let, for all $n\in\N$,
\[
\ki_n:=H(x)-H(z(n)), \nu_n:=\frac{J(z(n))}{H(x(n)) \ki_n},
\]
with the convention that $\nu_n:=0$ if $ \ki_n=0$.

By Lemma \ref{lem:hestimates} there exist $\epsilon$, $\lambda$,
$\mu :=\Cst(x,a)$ such that, for all $n\in\N$ such that $x(n)\in
B_{V_x}(2\epsilon)$, $\nu_n\in[\lambda,\mu]$. On the other hand,
for all $n\in\N$, using Lemma \ref{lem:stoc h} and the observation
that $J(z(n))=0$ if $ \ki_n=0$ by Lemma \ref{lem:hestimates},
\begin{eqnarray}\label{itka}
\ki_{n+1}&=& \biggl(1-\frac{\nu_n}{n+n_0+1} \biggr)
\ki_n-\xi_{n+1}-s_{n+1}\nonumber\\[-8pt]\\[-8pt]
&\le&\biggl(1-\frac{\lambda}{n+n_0+1} \biggr) \ki_n-\xi_{n+1}+s'_{n+1},
\nonumber
\end{eqnarray}
where
\[
s'_{n+1}:=-s_{n+1}+(\nu_n-\lambda)\max(- \ki_n,0)/(n+n_0+1).
\]
If $x(n)\in B_{V_x}(2\epsilon)$ for sufficiently small $\epsilon
:=\Cst(x,a)$, then, by Lemma \ref{lem:stoc h},
%
\begin{equation}
\label{ubd_xi}
\|\xi_{n+1}\|_\iy\le\frac{\Cst(x,a)}{n+n_0},\qquad  \|s'_{n+1}\|_\iy\le
\frac{\Cst(x,a)}{(n+n_0)^2},
\end{equation}
where we use in the second inequality that $\max(- \ki_n,0)\le\Cst
(x,a)/(n+n_0+1)$, since $\|x(n)-z(n)\|_\iy\le\Cst(x,a)/(n+n_0+1)$ by
Lemma \ref{bound_rn}(d), and $H(x(n))\le H(x)$ by Lemma \ref
{lem:hestimates}.

Let, for all $n\in\N$,
\[
\be_n:=\prod_{k=1}^n \biggl(1-\frac{\lambda}{k+n_0} \biggr).\vadjust{\goodbreak}
\]
Note that $\be_n n^\lambda$ converges to a positive limit. Inequality
(\ref{itka}) implies by induction that, for all $n\in\N$,
\[
\ki_n\le\be_n \biggl( \ki_0-\sum_{j=1}^n\frac{\xi_j}{\be_j}+\sum
_{j=1}^n\frac{s'_j}{\be_j} \biggr).
\]

Assume $\Lc(B_{V_x}(\epsilon))$ holds so that, in particular,
$x(n)\in\Lc(B_{V_x}(2\epsilon))$ for large $n\in\N$. The upper
bounds (\ref{ubd_xi}) yield, assuming w.l.o.g. $\lambda<1/2$, that
$\sum_{j=1}^n s'_j/\be_j<\iy$ and $\sum_{j=1}^n\Es(\xi_j^2)/\be
_j^2<\iy$; the latter implies, by the Doob convergence theorem in $\Lc
^2$, that $\sum_{j=1}^n \xi_j/\be_j$ converges a.s. Therefore, $ \ki
_n n^\lambda$ is bounded a.s.

We deduce subsequently, by Lemma \ref{lem:hestimates}(a), that
for all $\lambda\le\Cst(x,a)$,\break $J(x(n))n^\lambda$ converges a.s. to
$0$, so that $\lim_{n\to\iy}x(n)_{\partial S}n^\lambda=0$ in
particular. This implies that $\lim_{n\to\iy}I_{x(\iy
)}(x(n))n^\lambda=0$ by Lemma \ref{lem:hestimates}(b).

Now apply Lemma \ref{lem:evolent} with $q:=x(\iy)$: for large $n\in
\N$,
\begin{eqnarray*}
V_{x(\iy)}(z(n))&=&-\sum_{k=n}^\iy\frac{I_{x(\iy)}(x(k))}{k+n_0+1}+
\Biggl(x(\iy),\sum_{k=n+1}^\iy\ze_k \Biggr)
-2\sum_{k=n+1}^\iy(\epsilon_k)_{\partial S}\\
&&{} +\Cst(x,a)\Box\Biggl(\sum_{k=n}^\iy\frac{1}{(k+n_0)^2} \Biggr)\\
&=&o(n^{-\lambda})\qquad \mbox{a.s.},
\end{eqnarray*}
if we still assume w.l.o.g. $\lambda<1/2$, so that $\sum_{k=n+1}^\iy
(\epsilon_k)_{\partial S}=o(n^{-\lambda})$ a.s by Lemmas~\ref
{bound_rn}(a) and \ref{mart-est}. This completes the proof of
the lemma, using (\ref{eqVqiy1}).\vspace*{-3pt}

\subsection{\texorpdfstring{Proof of Lemma \protect\ref{lem:aborder}}{Proof of Lemma 10}}
\label{sec:aborder}
Let, for all $n\in\N$ and $i, j \in G$, $i\sim j$,
\[
Y_n^{i,j}:=\sum_{k=1}^n\frac{\1_{\{X_{k-1}=i,X_k=j\}}}{Z_{k-1}(j)},\qquad
Y_n^i:=\sum_{k=1}^n\frac{\1_{\{X_{k-1}=i\}}}{\sum_{j\sim
i}a_{j,i}Z_{k-1}(j)}.
\]
Then, by definition of the vertex-reinforced random walk,
\[
M_n^{i,j}:=Y_n^{i,j}-a_{i,j}Y_n^{i}
\]
is a martingale, and
%
\begin{eqnarray}\label{eq:ubdmij}
&&\sum_{k=1}^\iy\Es\bigl((M_k^{i,j}-M_{k-1}^{i,j})^2\bigr)\nonumber \\[-1pt]
&&\qquad =\Es\Biggl(\sum_{k=1}^\iy\frac{\1_{\{X_{k-1}=i\}}}{Z_{k-1}(j)^2}\frac
{a_{i,j}Z_{k-1}(j)}{\sum_{j\sim i}a_{j,i}Z_{k-1}(j)}
\biggl(1-\frac{a_{i,j}Z_{k-1}(j)}{\sum_{j\sim i}a_{j,i}Z_{k-1}(j)} \biggr)
\Biggr)\\[-1pt]
&&\qquad \le\Es\Biggl(\sum_{k=1}^\iy\frac{\1_{\{X_{k-1}=i,X_k=j\}
}}{Z_{k-1}(j)^2} \Biggr)<\iy\nonumber
\end{eqnarray}
so that, by the Doob convergence theorem in\vadjust{\goodbreak} $\Lc^2$, $M_n^{i,j}$
converges a.s.

Hence, for all $i\in\partial S$,
\begin{eqnarray*}
\log Z_n(i)&\equiv&\sum_{k=1}^n\frac{\1_{\{X_k=i\}
}}{Z_{k-1}(i)}=\sum_{j\sim i}Y_n^{j,i}\equiv\sum_{j\sim
i}a_{j,i}Y_n^{j}\\
&=&\sum_{j\sim i}a_{j,i}\sum_{k=1}^n\frac{\1_{\{X_{k-1}=j\}
}}{Z_{k-1}(j)}\frac{x(k-1)_j}{N_j(x(k-1))}\\
&\equiv&\sum_{j\sim i,j\notin\partial S}a_{i,j}\frac{x(\iy
)_j}{N_j(x(\iy))}\sum_{k=1}^n\frac{\1_{\{X_{k-1}=j\}}}{Z_{k-1}(j)}
\equiv\frac{N_i(x(\iy))}{H(x(\iy))}\log n,
\end{eqnarray*}
using Lemma \ref{lem:speed}, the symmetry of $a$ and $N_j(x(\iy))\not
=0$ for all $j\in G=T(x)$ in the third equivalence, and $H(x(\iy
))=N_j(x(\iy))$ for all $j\in S$ in the fourth equivalence [$x(\iy)$
being an equilibrium].

\subsection{\texorpdfstring{Proof of Proposition \protect\ref{prop:conv}}{Proof of Proposition 4}}
\label{sec:loc}
We will compare the probability of arbitrary paths remaining in $T(x)$
for the VRRWs defined, respectively, on the graphs~$T(x)$ and~$G$.
Let $x(n)$ [and its limit $x(\iy)$] denote the vector of occupation
density defined in the \hyperref[sec:gintro]{Introduction}, on the (finite)
subgraph~$T(x)$.\looseness=1

Let us introduce some notation. For all $k\in\N$ and $A\subset G$,
let $\Pc^A:=A^{\N}$ be the set of infinite sequences taking values in
$A$, and let $\T_k^A$ be the smallest $\s$-field on $\Pc^A$ that
contains the cylinders
\[
\C_{v,k}^A:=\{w\in\Pc^A\mbox{ s.t. } w_0=v_0,\ldots,w_k=v_k\},\qquad
v\in A^k.
\]
Let $\T^A:=\bigvee_{k\in\N}\T_k^A$. Finally, let $(X_j^A)_{j\ge n}$ be the VRRW on $A$ after time $n$,
conditionally to $X_n\in A$ (and be constant equal to $X_n$
otherwise).

For all $k\ge n$ and $v\in T(x)^k$,
\[
\Pb\bigl((X_{n+1},\ldots,X_k)=v | \F_n\bigr)=\Pb\bigl(\bigl(X_{n+1}^{T(x)},\ldots
,X_k^{T(x)}\bigr)=v | \F_n\bigr)Y_{n,k}^{(v)},
\]
where
%
\begin{equation}
\label{ynk}
Y_{n,k}:=\prod_{j=n}^{k-1}\prod_{\al\in\partial S(x)}
\biggl(1-\1_{\{X_j=\al\}}\frac{\sum_{\gamma\sim\al, \gamma\in
G\setminus T(x)} a_{\al,\gamma}Z_n(\gamma)}{\sum_{\be\sim\al
}a_{\al,\be}Z_j(\be)} \biggr)\in(0,1),\hspace*{-20pt}
\end{equation}
and $Y_{n,k}^{(v)}$ denotes the value of $Y_{n,k}$ at $(X_{n+1},\ldots
,X_k):=v$, where $Z_j(w)$, $w\in G$, $n\le j \le k-1$, assumes the
corresponding number of visits of $X_\cdot$ to~$w$.

We easily deduce that, for all $E\in\T^{T(x)}$,
\[
\Pb\bigl((X_{j+n})_{j\in\N}\in E | \F_n\bigr)=\Es\bigl(\1_{(X_{j+n}^{T(x)})_{j\in
\N}\in E}Y_{n,\iy} | \F_n\bigr).
\]

Let us now apply this equality with $E:=\{\Rc_{n,\iy}=T(x)\}\cap\Lc
(B_{V_x}(\epsilon))\cap\A_\partial(x(\iy))$ and prove that,
a.s.
on $E$, $Y_{n,\iy}>0$, which will complete the proof of the proposition:
for all $\al\in\partial S(x)$, a.s.\vadjust{\goodbreak} on $E$, if $\epsilon$ is
sufficiently small, then\looseness=1
\begin{eqnarray*}
\sum_{j=k}^\iy\frac{\1_{\{X_j=\al\}}}{\sum_{\be\sim\al}a_{\al
,\be}Z_j(\be)}
&=&\sum_{j=k}^\iy\frac{Z_j(\al)-Z_{j-1}(\al)}{\sum_{\be\sim\al
}a_{\al,\be}Z_j(\be)}\\[-1pt]
&\le&
\sum_{j=k}^\iy Z_j(\al) \biggl(\frac{1}{\sum_{\be\sim\al}a_{\al,\be
}Z_j(\be)}-\frac{1}{\sum_{\be\sim\al}a_{\al,\be}Z_{j+1}(\be)}
\biggr)\\[-1pt]
&\le&\oa\sum_{j=k}^\iy\frac{Z_{j}(\al)}{ (\sum_{\be\sim\al
}a_{\al,\be}Z_j(\be) )^2}\1_{\{X_{j+1}\sim\al\}}\\[-1pt]
&\le&\oa\sum_{j=k}^{\iy}\frac{x_j(\al)}{j(N_\al(x(j)))^2}<\iy,
\end{eqnarray*}\looseness=0
where we use that, since $\A_\partial(x(\iy))$ holds, $x(j)_\al\sim
_{j\to\iy}Cj^{N_\al(x(\iy))/H(x)-1}$ for some random $C>0$, so that
$\frac{x(j)_\al}{j(N_\al(x(j))^2}\sim_{j\to\iy}C\frac{j^{N_\al
(x(\iy))/H(x)-2}}{N_\al(x(\iy))}$, and $N_\al(x(\iy))<H(x(\iy
))=H(x)$ is $\epsilon$ is sufficiently small.\vspace*{-2pt}

\begin{appendix}
\section*{Appendix}
\subsection{Remainder of square-bounded martingales}
\label{sec:rsm}
The following lemma provides an almost sure estimate of $M_n-M_\iy$
for large $n$, when $M_n$ is a~martingale bounded in $L^2(\Om,\F,\Pb)$.
\begin{lem}
\label{mart-est}
Let $(M_n)_{n\ge0}$ be a bounded martingale in $L^2$, and let $f\dvtx\allowbreak \R
_+\rightarrow\R_+$ be a nondecreasing function such that $\int
_0^1(f(x))^{-2} \,dx<\iy$. Then
\[
M_n-M_\iy=o\bigl(f\bigl(\Es\bigl((M_n-M_\iy)^2\bigr)\bigr)\bigr)\qquad \mbox{a.s.}
\]
\end{lem}
\begin{pf}
For all $n\ge0$, let $s_n:=\Es((M_n-M_\iy)^2)$ and let
\[
N_n:=\sum_{k=1}^n\frac{M_k-M_{k-1}}{f(s_{k-1})},\qquad  N_0:=0.
\]
Then, for all $n\ge0$,
\[
\Es[N_n^2]=\sum_{k=1}^n\frac{s_{k-1}-s_k}{f(s_{k-1})^2}\le\int
_0^{s_0}\frac{dx}{(f(x))^2}<\iy.
\]
Therefore, $(M_n)_{n\ge0}$ and $(N_n)_{n\ge0}$ are martingales
bounded in $L^2$, and thus converge a.s.

Now, letting $O_n:=N_n-N_\iy$ for all $n\ge0$,
\begin{eqnarray*}
M_n-M_\iy&=&\sum_{k=n}^\iy f(s_k)(O_k-O_{k+1})
=f(s_n)O_n+\sum_{k=n+1}^\iy\bigl(f(s_k)-f(s_{k-1})\bigr)O_k\\
&=&o(f(s_n))\mbox{ a.s.}\vadjust{\goodbreak}
\end{eqnarray*}
\upqed
\end{pf}

\subsection{\texorpdfstring{Proof of Proposition \protect\ref{prop:cex}}{Proof of Proposition 1}}
\label{sec:cex}
Assume $X_0:=0$ for simplicity. Let, for all \mbox{$n\in\N$},
\begin{eqnarray*}
A_n&:=&Z_{n}(-1)+Z_{n}(1),\qquad\hspace*{12.5pt} \al_n^\pm:=Z_{n}(\pm1)/A_n,\\
R_n&:=&Z_n(0)/A_n-\log A_n, \qquad S_n:=\log\biggl(\frac{Z_{n}(-1)}{Z_{n}(1)}
\biggr)=\log\biggl(\frac{\al_n^-}{1-\al_n^-} \biggr).
\end{eqnarray*}

Let $a\in(0,1)$, $\epsilon<[a\wedge(1-a)]/2$. Given $n_0\in\N$
with $Z_{n_0}(0)$ sufficiently large and $X_{n_0}=0$, assume that
$Z_{{n_0}}(-2)/\log Z_{{n_0}}(-1)$, $Z_{{n_0}}(2)/\log Z_{{n_0}}(1) \in(1/3,1/2)$, $Z_{n_0}(\pm3)\le\Cst$, $\al^-_{n_0}\in
(a-\epsilon/3,a+\epsilon/3)$ and $R_{n_0}\in(-\epsilon/3,\epsilon
/3)$, which trivially occurs with positive probability.

Let us define the following stopping times:
\begin{eqnarray*}
T_0&:=&\inf\bigl\{n\ge{n_0}\mbox{ s.t. }X_n\in\{-3,3\}\mbox{ or
}X_n=X_{n-2}\in\{-2,2\}\bigr\},\\
T_1&:=&\inf\{n\ge{n_0}\mbox{ s.t. }Z_{n}(2)\vee Z_{n}(-2)>\log
Z_n(0)\},\\
T_2&:=&\inf\{n\ge n_0\mbox{ s.t. }\al^-_n\notin(a-\epsilon
/2,a+\epsilon/2)\mbox{ or }R_n\notin(-\epsilon/2,\epsilon/2)\},\\
T&:=&T_0\wedge T_1\wedge T_2.
\end{eqnarray*}

For all $n\in\N$, let $t_n$ be the $n$th return time to $0$, and let
$t'_n:=t_{n}\wedge T$.

As long as $n_0\le t_n<T$, $Z_{n}(0)=A_n(\log A_n+\Box(\epsilon/2))$,
which implies, for sufficiently large $Z_{n_0}(0)$, $A_{t_n}\le n$ by
contradiction, hence, $A_{t_n}\ge n/(\log n+\epsilon/2)$ and,
subsequently, $A_{t_n}\le n/(\log(\frac{n}{\log n+\epsilon
/2})-\epsilon/2)$.
Therefore, $Z_{t_n}(-1)\in((a-\epsilon)n/\log n,(a+\epsilon)n/\log
n)$ and $Z_{t_n}(1)\in((1-a-\epsilon)n/\log n,(1-a+\epsilon)n/\log
n)$, if $Z_{n_0}(0)\ge\Cst(a,\epsilon)$.

We successively upper bound $\Pb(T_0<T_1\wedge T_2 | \F_{{n_0}})$,
$\Pb(T_1<T_0\wedge T_2 | \F_{{n_0}})$ and $\Pb(T_2<T_0\wedge T_1 |
\F_{{n_0}})$, which will enable us to conclude that $\Pb(T=\iy|\F
_{n_0})>0$ for large $Z_{n_0}(0)$.

First, for sufficiently large $Z_{n_0}(0)$,
\begin{eqnarray}
\nonumber
&&\Pb(T_0<T_1\wedge T_2 | \F_{{n_0}})\\
\nonumber
&&\qquad \le\sum_{n\ge Z_{n_0}(0) : t_n<T}\Pb(X_{t_n+2}=X_{t_n+3}\mp
1=X_{t_n+4}=\pm2 | \F_{{n_0}})\nonumber\\
&&\qquad \quad {}+\Pb(X_{t_n+3}=\pm3 | \F_{{n_0}})\\
&&\qquad \le\Cst(a,\epsilon)\sum_{n\ge Z_{n_0}(0)} \biggl[\frac{1}{\log n}
\biggl(\frac{\log n}{n} \biggr)^2+
\frac{1}{\log n}\frac{\log n}{n}\frac{Z_{n_0}(3)+Z_{n_0}(-3)}{n/\log
n} \biggr]\nonumber \\
&&\qquad \le\Cst(a,\epsilon)\sum_{n\ge Z_{n_0}(0)}\frac{\log n}{n^2}<\frac{1}{3}.\nonumber
\end{eqnarray}

Let $\Gf:=(\F_{t'_n})_{n\ge Z_{n_0}(0)}$, and let us consider the
Doob decompositions of the $\Gf$-adapted processes $R_{t'_n}$ and
$S_{t'_n}$, $n\ge{Z_{n_0}(0)}$:
\[
R_{t'_n}=R_{n_0}+\De_n+\Psi_n,\qquad
S_{t'_n}:=S_{n_0}+\Phi_n+\Xi_n,
\]
where $\De_{Z_{n_0}(0)}=\Phi_{Z_{n_0}(0)}=\Psi_{Z_{n_0}(0)}=\Xi
_{Z_{n_0}(0)}:=0$ and, for all $n>{Z_{n_0}(0)}$,
\[
\De_n-\De_{n-1}:=\Es(R_{t'_n}-R_{t'_{n-1}} | \F_{t'_{n-1}}),\qquad  \Phi
_n-\Phi_{n-1}:=
\Es(S_{t'_n}-S_{t'_{n-1}} | \F_{t'_{n-1}}),
\]
and $(\Psi_n)_{n\ge Z_{n_0}(0)}$ and $(\Xi_n)_{n\ge Z_{n_0}(0)}$ are
$\Gf$-adapted martingales.

Let us now estimate the expectation and variance of the increments of
the processes $(R_{t'_n})_{n\in\N}$: if $n\ge\Cst(\epsilon)$,
\begin{eqnarray*}
\Es(R_{t'_n+1}-R_{t'_n} | \F_{t'_{n}})
&=&\frac{1}{A_{t'_n}}+\frac{A_{t'_n}}{A_{t'_n}+n} \biggl(-\frac
{n+1}{A_{t'_n}(A_{t'_n}+1)}-\frac{1}{A_{t'_n}}+\Box\biggl(\frac
{1}{(A_{t'_n})^2} \biggr) \biggr)\\
&=&\frac{1}{A_{t'_n}}+\frac{A_{t'_n}}{A_{t'_n}+n} \biggl(-\frac
{n}{A_{t'_n}(A_{t'_n}+1)}-\frac{1}{A_{t'_n}} \biggr)+\Box\biggl(\Cst\frac{\log
n}{n^2} \biggr)\\
&=&-\frac{n}{A_{t'_n}(A_{t'_n}+1)(A_{t'_n}+n)}+\Box\biggl(\Cst\frac{\log
n}{n^2} \biggr)\\
&=&\Box\biggl(\Cst\frac{(\log n)^3}{n^2} \biggr),
\\
\Es(R_{t'_{n+1}}-R_{t'_n+1} | \F_{t'_{n}})&=&\Box\biggl(\Cst\frac{\log
n}{n} \biggl(\frac{1}{A_{t'_n}}+\frac{n+1}{(A_{t'_n}+1)(A_{t'_n}+2)} \biggr)
\biggr)\\
&=&\Box\biggl(\Cst\frac{(\log n)^3}{n^2} \biggr)
\end{eqnarray*}
and
\[
|R_{t'_{n+1}}-R_{t'_n}|\le\frac{2}{A_{t'_n}}+\frac
{2(n+1)}{A_{t'_n}(A_{t'_n}+1)}\le\Cst\frac{(\log n)^2}{n}
\]
so that, in summary,
%
\begin{equation}
\label{eq:pbd1}
\qquad |\De_n-\De_{n-1}|\le\Cst\frac{(\log n)^3}{n^2},\qquad  \Es\bigl((\Psi
_{n+1}-\Psi_n)^2 | \F_{t'_{n}}\bigr)\le\Cst
\frac{(\log n)^4}{n^2}.
\end{equation}
Let us do similar computations for $(S_{t'_n})_{n\in\N}$: if $n\ge
\Cst(a,\epsilon)$,
\begin{eqnarray*}
\Phi_n-\Phi_{n-1}&=& \biggl(\frac{1}{Z_{t_n}(-1)}+\frac
{1}{Z_{t_n}(-1)}\Box\biggl(\frac{\log n}{n} \biggr)+\Box\biggl(\frac
{1}{(Z_{t_n}(-1))^2} \biggr) \biggr)\\
&&{}\times\frac{Z_{t_n}(-1)}{n+Z_{t_n}(-1)+Z_{t_n}(1)}\\
&&{} - \biggl(\frac{1}{Z_{t_n}(1)}+\frac{1}{Z_{t_n}(1)}\Box\biggl(\frac
{\log n}{n} \biggr)+\Box\biggl(\frac{1}{(Z_{t_n}(1))^2} \biggr) \biggr)\\
&&\quad {}\times
\frac{Z_{t_n}(1)}{n+Z_{t_n}(-1)+Z_{t_n}(1)}\\
&=&\frac{\log n}{n^2}\Box(\Cst(a,\epsilon) ),
\end{eqnarray*}
and
\[
|S_{t'_{n+1}}-S_{t'_n}|\le\log\biggl(1+\frac{2}{Z_{t_n}(1)} \biggr)\vee\log
\biggl(1+\frac{2}{Z_{t_n}(-1)} \biggr),
\]
so that
%
\begin{eqnarray}
\label{eq:pbd2}
|\Phi_n-\Phi_{n-1}|&\le&\Cst(a,\epsilon)\frac{\log n}{n^2},\qquad\\
\Es\bigl((\Xi_{n+1}-\Xi_n)^2 | \F_{t'_{n}}\bigr)&\le&\Cst(a,\epsilon) \biggl(\frac
{\log n}{n} \biggr)^2.\nonumber
\end{eqnarray}
Hence, by Chebyshev's and Doob's martingale inequalities, for all
$\delta>0$,
\[
\Pb\Bigl(\max_{k\ge Z_{n_0}(0)}|\Psi_k|>\delta| \F_{n_0} \Bigr)\le\frac
{\Cst}{\delta^2}\sum_{j=Z_{n_0}(0)}^\iy
\frac{(\log n)^4}{n^2}
\le\frac{\Cst}{\delta^2}\frac{(\log Z_{n_0}(0))^4}{Z_{n_0}(0)}
\]
and a similar inequality holds on the maximum of $|\Xi_k|$, $k\ge
Z_{n_0}(0)$, so that, for sufficiently large $Z_{n_0}(0)$, $\Pb
(T_2<T_0\wedge T_1 | \F_{{n_0}})<1/3$.

Let us now make use of notation $Y_n^{i,j}$, $Y_n^i$ and $M_n^{i,j}$
from Section \ref{sec:aborder} (with $a_{i,j}=\1_{i\sim j}$), and
let $U_n^\pm:=Y_n^{\pm1,\pm2}$, $V_n^\pm:=Y_n^{\pm1}$ and $W_n^\pm
:=M_n^{\pm1,\pm2}=U_n^\pm-V_n^\pm$. Then the processes $(U_n^\pm
)_{n\ge0}$ are martingales and, using \eqref{eq:ubdmij}, for all
$n\ge n_0$,
%
\begin{eqnarray}
\label{eq:pbd3}
\Es\bigl((W_n^\pm-W_{n_0}^\pm)^2 | \F_{{n_0}}\bigr)&\le&
\Es\Biggl(\sum_{k=n_0+1}^n\frac{\1_{\{X_{k-1}=\pm1,X_k=\pm2\}}}{Z_{k-1}(\pm2)^2}\Big|\F_{{n_0}}\Biggr)
\nonumber\\[-8pt]\\[-8pt]
&\le&\sum_{j\ge Z_{n_0}(\pm2)}\frac{1}{j^2}\nonumber
\end{eqnarray}
so that, if $\Upsilon:=\{\max_{k\ge n_0}|W_k^i-W_{n_0}^i|\le\delta
,i\in\{+,-\}\}$, then, for all $\delta>0$,
\[
\Pb(\Upsilon^c | \F_{n_0} )\le\frac{1}{\delta^2} \biggl(\frac
{1}{Z_{n_0}(2)-1}+\frac{1}{Z_{n_0}(-2)-1} \biggr)<\frac{1}{3}
\]
for sufficiently large $Z_{n_0}(0)$.

Now, on $\Upsilon$, for all $n<T$, choosing $\delta=(\log2)/3$, and
again for sufficiently large $Z_{n_0}(0)$,
\begin{eqnarray*}
\log Z_n(\pm2)&\le&\log Z_{n_0}(\pm2)+U_n^\pm- U_{n_0}^\pm+\delta
\le2\delta+\log Z_{n_0}(\pm2)+V_n^\pm-V_{n_0}^\pm\\
&\le&2\delta+\log Z_{n_0}(\pm2)+\sum_{k=n_0+1}^n\frac{\1_{\{
X_{k-1}=\pm1\}}}{Z_{k-1}(0)}\\
&\le&2\delta+\log Z_{n_0}(\pm2)+\sum_{k=Z_{n_0}(\pm1)}^{Z_{n-1}(\pm
1)}\frac{1}{k\log k}\\
&\le&3\delta+\log\biggl(\frac{Z_{n_0}(\pm2)}{\log Z_{n_0}(\pm1)} \biggr)+\log
(\log Z_n(\pm1))\\
&\le&\log(\log Z_n(\pm1))\le\log(\log Z_n(0)),
\end{eqnarray*}
where we use in the fourth inequality that, if $n<T$, then $T_n\ge
-\epsilon/2$ and $\al^-_n\in(a-\epsilon/2,a+\epsilon/2)$ so that
$Z_n(0)\ge Z_n(\pm1)\log Z_n(\pm1)$ if $Z_{n_0}\ge\Cst(a,\epsilon
)$, and in the sixth inequality that $Z_{n_0}(\pm2)/\log Z_{n_0}(\pm
1)\le1/2$. This completes the proof, as $\Pb(T_1<T_0\wedge T_2 | \F
_{{n_0}})\le\Pb(\Upsilon^c|\F_{n_0})<1/3$ for large $Z_{n_0}(0)$.

The estimates \eqref{eq:pbd1}--\eqref{eq:pbd2} [resp., \eqref
{eq:pbd3}] imply that the $\Gf$ (resp., $\Ff$)-adapted martingales
$(\Psi_n)_{n\ge Z_{n_0}(0)}$ and $(\Xi_n)_{n\ge Z_{n_0}(0)}$ (resp.,
$W_n^\pm$) are bounded in $L^2$ and hence converge a.s.

Therefore, on $\{T=\iy\}$, (i)--(ii) hold, and $(\al_n)_{n\ge0}$ and
$(R_n)_{n\ge0}$ converge a.s. Note that Lemma \ref{mart-est} implies
more precisely, for all $\nu<1/2$, $\Xi_n-\Xi_\iy=o(n^{-\nu})$,
hence, $\al_n-\al_\iy=o(Z_n(0)^{-\nu})$. Thus, on $\{T=\iy\}$,
\begin{eqnarray*}
\log Z_n(\pm2)&\equiv& U_n^\pm\equiv V_n^\pm=\sum_{k=0}^{n-1}\frac
{\1_{\{X_k=\pm1\}}}{Z_k(\pm2)+Z_k(0)}\\
&\equiv&\al_\iy^\pm\sum_{k=0}^{n-1}\frac{\1_{\{X_k=\pm1\}
}}{Z_k(\pm1)\log Z_k(\pm1)} \biggl(1+O \biggl(\frac{1}{\log Z_k(\pm1)} \biggr) \biggr)\\
&\equiv&\al_\iy^\pm\log(\log Z_n(\pm1))\equiv\al_\iy^\pm\log
(\log n),
\end{eqnarray*}
which proves (iii).
\end{appendix}

\section*{Acknowledgments}

We would like to thank the referees for very helpful comments.


%

\printaddresses

\end{document}